\documentclass[11pt,a4paper]{article}%
\usepackage{amsmath,amssymb}
\usepackage{xcolor}
\usepackage{amsmath}%
\usepackage{amsfonts}%
\usepackage{amssymb}%
\usepackage{graphicx}
\providecommand{\U}[1]{\protect\rule{.1in}{.1in}}
\newcommand{\p}{\partial}

\newcommand{\R}{\mathbb{R}}

\def\tilde{\widetilde}

\newtheorem{theorem}{Theorem}
\newtheorem{proposition}{Proposition}

\newtheorem{notation}{Notation}
\newtheorem{remarka}{Remark}

\newtheorem{lemme}{Lemma}

\allowdisplaybreaks

\begin{document}

\title{New effective pressure and existence of global strong solution for
compressible Navier-Stokes equations with general viscosity coefficient in one dimension}
\author{Cosmin Burtea \thanks{Universit\'e Paris Diderot, Sorbonne Paris Cit\'e, Institut de Math\'ematiques de Jussieu-Paris Rive Gauche (UMR 7586), F-75205 Paris, France}, Boris Haspot \thanks{Universit\'e Paris Dauphine, PSL
Research University, Ceremade, Umr Cnrs 7534, Place du Mar\' echal De Lattre
De Tassigny 75775 Paris cedex 16 (France), haspot@ceremade.dauphine.fr }
\thanks{ANGE project-team (Inria, Cerema, UPMC, CNRS), 2 rue Simone Iff, CS
42112, 75589 Paris, France.}}
\date{}
\maketitle

\begin{abstract}
In this paper we prove the existence of global strong solution for the
Navier-Stokes equations with general degenerate viscosity coefficients. The
cornerstone of the proof is the introduction of a new effective pressure which
allows to obtain an Oleinik-type estimate for the so called effective
velocity. In our proof we make use of additional regularizing effects on the
velocity which requires to extend the technics developed by Hoff for the
constant viscosity case.

\end{abstract}


\section{Introduction}

We consider the compressible Navier Stokes system in one dimension with
$x\in\mathbb{R}$:
\begin{equation}%
\begin{cases}
\begin{aligned} &\partial_{t}\rho+\partial_x(\rho u)=0,\\ &\partial_{t}(\rho u)+\partial_x(\rho u^2)-\partial_x(\mu(\rho)\partial_x u)+\partial_x P(\rho)=0,\\ &(\rho,u)_{t=0}=(\rho_0,u_0). \end{aligned}
\end{cases}
\label{0.1}%
\end{equation}
Here $u=u(t,x)\in\mathbb{R}$ stands for the velocity field, $\rho=\rho
(t,x)\in\mathbb{R}^{+}$ is the density, $P(\rho)=\rho^{\gamma}$ is the
pressure.
We denote by $\mu(\rho)$ the viscosity coefficient of the fluid and $(\rho
_{0},u_{0})$ are the initial data. In the sequel we shall only consider
viscosity of the form:
\begin{equation}
\mu(\rho)=\rho^{\alpha}\label{condition}%
\end{equation}
with $\alpha>0$. This choice is motivated by
physical considerations. Indeed it is justified by the derivation of the
Navier-Stokes equations from the Boltzmann equation through the Chapman-Enskog
expansion to the second order (see \cite{CC70}), the viscosity coefficient is
then a function of the temperature. If we consider the case of isentropic
fluids, this dependence is expressed by a dependence on the density function (we refer in particular to \cite{HoffSerre}).
We mention that the case $\mu(\rho)=\rho$ is
related to the so called viscous shallow water system. This system with
friction has been derived by Gerbeau and Perthame in \cite{GP} from the
Navier-Stokes system with a free moving boundary in the shallow water regime
at the first order. This derivation relies on the hydrostatic approximation
where the authors follow the role of viscosity and friction on the bottom.
\newline We are now going to rewrite the system (\ref{0.1}) following the new
formulation proposed in \cite{para} (see also \cite{CPAM,CPAM1,MN}), indeed
setting:
\begin{equation}
v=u+\frac{\mu(\rho)}{\rho^{2}}\partial_{x}\rho,\label{effect}%
\end{equation}
called the effective velocity, we can rewrite the system (\ref{0.1}) as
follows:
\begin{equation}%
\begin{cases}
\begin{aligned} &\partial_t\rho+\partial_x(\rho u)=0,\\ &\rho\partial_t v+\rho u\partial_x v+\partial_x P(\rho)=0. \end{aligned}
\end{cases}
\label{0.1a}%
\end{equation}
\indent The existence of global \textit{weak} solution has been obtained by Jiu and
Xin in \cite{Jiu} for viscosity coefficients verifying (\ref{condition}). In
passing we point out that a large amount of literature is essentially
dedicated to the study of the compressible Navier-Stokes equations with
constant viscosity coefficients. In particular the existence of global strong
solution with large initial data for initial density far away from the vacuum
has been proved for the first time by Kanel \cite{Ka} (see also
\cite{KS77,Hof87} ). In \cite{HS01} the authors proved that vacuum states do
not arise provided that the initial density is positive almost everywhere. We
would like also to mention the results of Hoff in \cite{Hof98} who proved the
existence of global weak solution for constant viscosity coefficients with
initial density admitting shocks (we refer also to \cite{Ser86a,Ser86b,Lio98}%
). The author exhibited regularizing effects on the velocity via the use of
tricky estimates on the convective derivative:
\[
\dot{u}=\partial_{t}u+u\partial_{x}u,
\]
we will generalize these techniques in the present paper to the case of general
viscosity coefficients. In \cite{mesure}, the second author proved also the
existence of global weak solution for general viscosity coefficients with
initial density admitting shocks and with initial velocity belonging to the
set of finite measures. In opposite to \cite{Hof98}, the initial data satisfy
the BD entropy but not the classical energy, it allows in particular to show
some regularizing effects on the density inasmuch as the density becomes
instantaneously continuous. It is due to the regularity of the effective
velocity $v$ which express the coupling between the velocity and the
density.
\newline\indent The problem of existence of global strong solution for system
(\ref{0.1}) with large initial data and with general viscosity coefficients
verifying (\ref{condition}) is not yet completely solved. Indeed when
$\alpha>1$ it requires conditions of sign on the so called effective flux (see
\cite{Hof87,Lio98}). This quantity represents the force that
the fluids exerts on itself and a priori has no reason to be signed. In the
following we are going to present the current state of art concerning the
existence of global strong solution for system (\ref{0.1}) with viscosity
coefficients verifying (\ref{condition}).
\newline\indent It has been first proved by
Mellet and Vasseur (see \cite{MV06}) in the case $0<\alpha<\frac{1}{2}$. The
main argument of their proof consists in using the Bresch-Dejardins  entropy (see \cite{BD})
in order to estimate the $L^{\infty}$ norm of $\frac{1}{\rho}$ and using the
parabolicity of the momentum equation of (\ref{0.1}). It is important at this
level to point out that the Bresch-Dejardins entropy gives almost for free the control of
$\Vert\frac{1}{\rho}\Vert_{L_{t,x}^{\infty}}$ when $\alpha<\frac{1}{2}%
$.
\newline \indent In \cite{MN}, the second author has proved similar results for the
case $\frac{1}{2}<\alpha\leq1$ where he exploited the fact that the effective
velocity $v$ satisfies a damped transport equation. It enables to obtain
$L^{\infty}$ estimates for $v$ and using maximum principle to get $L^{\infty}$
control on $\frac{1}{\rho}$.
\newline More recently Constantin \textit{et al} in \cite{Constantin} have extended the previous results. More precisely, in
the range $\alpha\in\left(  \frac{1}{2},1\right]  $ under the condition
$\gamma\geq2\alpha$, the authors obtain global existence of strong solutions
for initial data belonging to $H^{3}$. They prove that the same result also
holds true in the case $\alpha>1$ with $\gamma$ belonging to $[\alpha
,\alpha+1]$ provided that the initial data satisfy:
\begin{equation}
\partial_{x}u_{0}\leq\rho_{0}^{\gamma-\alpha}.\label{flux5}%
\end{equation}
We point out that the condition (\ref{flux5}) is equivalent to consider a
negative effective flux (see for example \cite{Hof87,Lio98}) at initial
time. The main idea of their proof consists in proving via a maximum principle
that the effective flux remains negative for all time. This is sufficient to
control the $L^{\infty}$ norm of $\frac{1}{\rho}$.
\newline \indent In the present paper, our goal is double inasmuch as we wish both to
show the existence of global strong solution for the case $\alpha>\frac{1}{2}$
without any sign restriction on the initial data and with minimal assumptions
in terms of regularity. In \cite{Constantin}, Constantin \textit{et al} proved a
blow-up criterion for $\alpha>\frac{1}{2}$ which is relied to estimating the
$L_{t,x}^{\infty}$ norm of $\frac{1}{\rho}$. In order to apply this blow-up
criterion, we introduce a new effective pressure $y=\frac{\partial_{x}v}{\rho
}+F_{2}(\rho)$ with $\rho F_{2}^{\prime}(\rho)=\frac{F_{1}(\rho)}{\rho}$ and
$F_{1}(\rho)=\frac{P^{\prime}(\rho)\rho}{\mu(\rho)}$. We observe then that $y$
satisfies the following equation:
\begin{equation}
\begin{aligned} &\partial_t y+u\partial_x y+F_1(\rho)y-F_1(\rho)F_2(\rho)+F_1' (\rho)\frac{\rho}{\mu(\rho)}(v-u)^2=0. \end{aligned}\label{gcru}%
\end{equation}
This last equation enables us to prove that if $y_{0}\leq C$ with
$C\in\mathbb{R}$ then $y$ remains bounded on the right all along the time
which implies in particular that:
\begin{equation}
\partial_{x}v(t,x)\leq C_{1}(t)\;\;\;\forall(t,x)\in\mathbb{R}^{+}%
\times\mathbb{R},\label{Oleinik1}%
\end{equation}
with $C_{1}$ a continuous increasing function. Using maximum principle for the
mass equation of (\ref{0.1a}) allows us to prove that $\frac{1}{\rho}$ is
bounded all along the time. In order to show the uniqueness of the solutions,
we extend Hoff's techniques to the case of general viscosity coefficients which
enables us to prove that $\partial_{x}u$ belongs to $L_{loc}^{1}(L^{\infty
}(\mathbb{R}))$. Passing in Lagrangian formulation (see the Appendix and the
references therein), we get the uniqueness of the solutions. Finally, we would
like to mention that the estimate (\ref{Oleinik1}) is reminiscent of the
so-called Oleinik estimate (see \cite{Oleinik,Daf}) for scalar
conservation law with a flux strictly convex or concave. If we consider the
following equation with $f$ regular:
\[
\partial_{t}u+\partial_{x}f(u)=0,\;u(0,\cdot)=u_{0}\in L^{\infty}(\mathbb{R}),
\]
the Kruzhkov theorem (see \cite{Kruz}) asserts that there exists a unique entropy
solution. In addition if $f$ is genuinely non linear, Oleinik has proved the
following estimate in the sense of measures for $C>0$ and for any
$(t,x)\in\mathbb{R}^{+}\times\mathbb{R}$:
\begin{equation}
\partial_{x}u(t,x)\leq\frac{C}{t}.\label{Oleinik2}%
\end{equation}
This estimate gives regularizing effects on $u$ since instantaneously
$u(t,\cdot)$ with $t>0$ is in $BV_{loc}(\mathbb{R})$. In our case, we have no
regularizing effects on $v$. A possible explanation is the fact that $v$
satisfies a damped transport equation which is in some sense linearly
degenerate.

\section{Main results}

We are now in position to state our main theorem.

\begin{theorem}
\label{theo1} Let $\alpha>\frac{1}{2}$, $\gamma\geq\max(1,\alpha)$, $(\rho_{0},\frac
{1}{\rho_{0}})\in(L^{\infty}(\mathbb{R}))^{2}$, $(\rho_{0}-1,u_{0})\in
(L^{2}(\mathbb{R}^{2}))^{2}$. In addition we assume that $v_{0}\in
L^{2}(\mathbb{R})$ and that there exists
$C\in\mathbb{R}$ such that for any $x>y$ we have:
\begin{equation}
\frac{v_{0}(x)-v_{0}(y)}{x-y}\leq C\label{BV}%
\end{equation}
Then there exists a unique global strong solution $(\rho,u)$ for the Navier-Stokes system $\eqref{0.1}$ with the
following properties. For any given $T>0$, $L>0$ there exist a positive
constant $C(T)$, a positive constant $C(T,L)$ depending respectively on $T$, on
$T$, $L$ and on $\Vert\rho_{0}-1\Vert_{L^{2}}$, $\Vert(\rho_{0},\frac{1}%
{\rho_{0}})\Vert_{L^{\infty}}$, $\Vert u_{0}\Vert_{L^{2}}$, $\Vert v_{0}%
\Vert_{L^{2}}$ such that, if $\sigma(t)=\min(1;t)$, then:
\begin{equation}
C(T)^{-1}\leq\rho(T,\cdot)\leq C(T)\;\;\;\mbox{a.e},\label{1.17}%
\end{equation}%
\begin{equation}
\begin{aligned} &\sup_{0<t\leq T}\big( \|\rho(t,\cdot)-1\|_{L^2}+\|u(t,\cdot)\|_{L^2}+\|\partial_x\rho(t,\cdot)\|_{L^2}+\sigma(t)^{\frac{1}{2}}\|\partial_x u(t,\cdot)\|_{L^2}\\ &+\sigma(t)^{\frac{1}{2}}(\|\dot{u}(t,\cdot)\|_{L^2}+\|\partial_x(\rho^\alpha\partial_x u(t,\cdot)-P(\rho)+P(1))\|_{L^2}\big)\leq C(T), \end{aligned}\label{1.18}%
\end{equation}%
\begin{equation}
\begin{aligned} &\int^T_0 [\|\partial_x u(t,\cdot)\|_{L^2}^2+\|\partial_x\rho(t,\cdot)\|_{L^2}^2+\sigma(t)\|\dot{u}(t,\cdot)\|_{L^2}^2+\sigma(t)\|\partial_x \dot{u}(t,\cdot)\|_{L^2}^2]dt \leq C(T), \end{aligned}\label{1.19}%
\end{equation}%
\begin{equation}
\int_{0}^{T}\sigma^{\frac{1}{2}}\left(  \tau\right)  \left\Vert \partial
_{x}u\left(  \tau\right)  \right\Vert _{L^{\infty}}^{2}d\tau\leq C\left(
T\right)  .\label{Lipo}%
\end{equation}%
\begin{equation}
\sup_{0<t\leq T}\sigma(t)^{\frac{1}{2}}\Vert\partial_{x}u(t,\cdot
)\Vert_{L^{\infty}}\leq C\left(  T\right)  \label{Lipo1}%
\end{equation}%
\begin{equation}
\Vert v\Vert_{BV([0,T]\times\lbrack-L,L])}\leq C(T,L).\label{BV2}%
\end{equation}
Furthermore for any $x>y$ and $t\geq 0$, we have almost everywhere:
\begin{equation}
\frac{v(t,x)-v(t,y)}{x-y}\leq C_1(t),\label{uBV}%
\end{equation} 
with $C_1$ a continuous increasing function.
\end{theorem}

\begin{remarka}
It is important to point out that our theorem requires that $\p_x\rho_0$ belongs to $L^2(\R)$.
Indeed since $\partial_{x} \varphi(\rho_{0})=v_{0}-u_{0}$ with $v_{0}\in
L^{2}(\mathbb{R})$ and $u_{0}\in L^{2}(\mathbb{R})$, it implies that
$\partial_{x} \varphi(\rho_{0})\in L^{2}(\mathbb{R})$.  since $\frac{1}{\rho_{0}}$ is in $L^{\infty}(\mathbb{R})$, we deduce that $\p_x\rho_0$ is in $L^2(\R)$.\\
Furthermore since $\varphi(\rho
_{0})-1$ is also in $L^{2}(\mathbb{R})$ using that $\frac{1}{\rho_{0}}$ and $\rho_{0}$ are in $L^{\infty}(\mathbb{R})$, we deduce that $\varphi(\rho_{0})-1$
is in $H^{1}(\mathbb{R})$. The initial density $\rho_{0}$ is then necessary a continuous
function which prevents us from considering shock-type initial data. 

\end{remarka}
\begin{remarka}
We would like to mention that any solution $(\rho,u)$ of system (\ref{0.1}) in the sense of distributions which verifies the regularity assumptions of Theorem \ref{theo1} is also a strong solution i.e. $(\rho,u)$ satisfy the system (\ref{0.1}) almost everywhere on $\R^+\times\R$.
Setting $w_1(t,x)=\rho^\alpha\partial_x u(t,x)-P(\rho(t,x))+P(1)$ the effective flux, we get from (\ref{1.17}) and (\ref{1.18}) that for any $t>0$:
$$
\begin{cases}
\begin{aligned}
&\sigma(t)^{\frac{1}{2}} \|\partial_x w_1(t,\cdot)\|_{L^2}\leq C(t)\\
&\sigma(t)^{\frac{1}{2}}\|\partial_x u(t,\cdot)\|_{L^2}+\|P(\rho)(t,\cdot)-P(1)\|_{L^2}\leq C(t),
\end{aligned}
\end{cases}
$$
for $C$ a continuous increasing function. This implies that $w_1$ belongs to $L^1_{loc}(\R^+,H^1(\R))$. Using now the fact that $(P(\rho)-P(1))$ belongs to $L^\infty_{loc}(H^1(\R))$, we deduce that $\rho^\alpha \p_x u$ is in $L^1_{loc}(H^1(\R))$. Using (\ref{1.17}), the fact that $(\frac{1}{\rho^\alpha}-1)$ belongs to $L^\infty_{loc}(H^1(\R))$ we get using product law in Sobolev spaces that $\p_x u$ is in $L^1_{loc}(H^1(\R))$. In particular $\p_{xx}u$ is in $L^1_{loc}(L^2(\R))$. In other words it is easy to observe that each term of (\ref{0.1}) is in $L^1_{loc}(\R^+\times\R)$ which ensures that $(\rho,u)$ satisfies (\ref{0.1}) almost everywhere.
\end{remarka}

\begin{remarka}
Let us point out that compared with \cite{Constantin}, we deal with the range
$\gamma\geq\max(\alpha,1)$ whereas in \cite{Constantin} the authors treat the case
$\alpha\leq\gamma\leq\alpha+1$, $\alpha>1$ provided that $\partial_{x}%
u_{0}\leq\rho_{0}^{\gamma-\alpha}$. In a certain sense the method that we developed in our proof unifies the different situations, $\gamma>\alpha+1$ and $\alpha\leq\gamma<\alpha+1$. Furthermore we do not require any condition of sign on the initial data.
\end{remarka}

\begin{remarka}
\label{remimpo}
The condition (\ref{BV}) is a condition of Oleinik-type which implies that
$v_{0}$ is in $BV_{loc}(\mathbb{R})$. Indeed we recall that for any
$x\in\mathbb{R}$ we have $|x|=(2x)_{+}-x$ with $(x)_{+}=\max(0,x)$. It yields
then that for any interval $[a,b]$ such that $v_{0}(a)$ and $v_{0}(b)$ are
finite and any increasing subdivision $(x_{n})_{n=1,\cdots,N}$ of the interval
$[a,b]$ with $N\in\mathbb{N}^{*}$ , we have using (\ref{BV}) and taking
$x_{0}=a$, $x_{N+1}=b$ if $x_{1}>a$ and $x_{N}<b$:
\[
\begin{aligned}
&\sum_{i=1}^{N-1}|v_0(x_{i+1})-v_0(x_i)|\leq \sum_{i=0}^{N}|v_0(x_{i+1})-v_0(x_i)|\\
&\leq 2\sum_{i=0}^{N}(v_0(x_{i+1})-v_0(x_i))_+ +v_0(a)-v_0(b)\\
&\leq 2C\sum_{i=0}^{N}(x_{i+1}-x_i) +v_0(a)-v_0(b)\\
&\leq 2C(b-a)+v_0(a)-v_0(b)\\
\end{aligned}
\]
In particular this shows that $v_0$ is necessary in $L^\infty_{loc}(\R)$.\\
Furthermore (\ref{uBV}) implies that the Oleinik estimate (\ref{BV}) is preserved all along the time.
In addition since $x\rightarrow v(t,x)-C_1(t)x$ is non-increasing, we deduce that $v(t,\cdot)$ has left and right-hand limits at each points for almost $t\geq 0$.
\end{remarka}

\begin{remarka}
Our theorem does not require high regularity assumption on the initial velocity. Indeed, we assume only that $u_0$ and $v_0$ are respectively in $L^2(\R)$ and $L^2(\R)\cap BV_{loc}(\R)$. This is however sufficient in order to ensure uniqueness.
\end{remarka}

\begin{remarka}
We can observe that in the case $\frac{1}{2}<\alpha\leq 1$, our assumption on $\gamma$ is optimal (from an hyperbolic point of view) since we need only $\gamma\geq 1$. This extends the results of \cite{Constantin,MN}.
\end{remarka}

\begin{remarka}
We can observe that (\ref{Lipo}) and (\ref{Lipo1}) give a $L^{1}%
_{loc}(L^{\infty}(\mathbb{R}))$ control on $\partial_{x} u$. In
particular, this enables us to define the flow associated to the velocity $u$ (we refer for
more details to the Appendix).
\end{remarka}
We would like to emphasize that the condition $\left( \text{\ref{BV}}\right) $ is
automatically satisfied provided that $\partial _{x}v_{0}\in L^{\infty }$.
A necessary condition for this later condition to hold is to take initial data $\left( \dfrac{1}{\rho _{0}}-1,\rho _{0}-1,u_{0}\right)$ in the following Sobolev spaces $\left( H^{s}\left( \mathbb{R}\right) \right)^{2}\times  H^{s-1}\left( \mathbb{R}\right) $  with $s_1>\frac{5}{2}$, . As a by-product of Theorem \ref{theo1} and the Appendix, we establish the following result.
\begin{theorem}
Consider $\alpha \geq \frac{1}{2}$, $\gamma \geq \max \left( 1,\alpha
\right) $ and  
$$
\left( \dfrac{1}{\rho _{0}}-1,\rho _{0}-1,u_{0}\right) \in\left( H^{s}\left( \mathbb{R}\right) \right)^{2}\times  H^{s-1}\left( \mathbb{R}\right) 
$$
with $s>\frac{5}{2}$.
Then, the compressible Navier-Stokes system $\left( \text{\ref{0.1}}\right) $ admits an
unique solution,  we have 
$$\left( \rho -1,u\right) \in C(\mathbb{R}_{+},H^{s}\left(
\mathbb{R}\right) \times H^{s-1}\left( 
\mathbb{R}\right)  ).$$
\end{theorem}

In the section \ref{section3}, we prove the Theorem \ref{theo1}. An appendix
is devoted to the definition and basic properties of the Lagrangian framework, we give also a sketch of the proof of the Theorem \ref{Cons} below.

\section{Proof of the Theorem \ref{theo1}}

\label{section3} 

A first ingredient is the following blow-up criterion 

\begin{theorem}
\label{Cons} Assume that $\alpha>\frac{1}{2}$ and $\gamma\geq\max(\alpha
-\frac{1}{2},1)$ and let $s\geq 3$ and $(\rho_{0}-1,u_{0})\in H^{s}(\mathbb{R})$.
Then there exists $T^{*}>0$ such that $(\rho,u)$ is a strong solution on
$(0,T^{*})$ with:
\[
(\rho-1)\in C(0,T,H^{s}(\mathbb{R})),\,u\in C(0,T,H^{s}(\mathbb{R}))\cap
L^{2}(0,T,H^{s+1}(\mathbb{R})),\;\forall T\in(0,T^{*}),
\]
and for all $t\in(0,T^{*})$:
\[
\|\frac{1}{\rho}(t,\cdot)\|_{L^{\infty}}\leq C(t),
\]
where $C(t)<+\infty$ if $t\in(0,T^{*})$. In addition, if:
\[
\sup_{t\in(0,T^{*})}\|\frac{1}{\rho}(t,\cdot)\|_{L^{\infty}}\leq C<+\infty,
\]
then the solution can be continued beyond $(0,T^{*})$.
\end{theorem}

The above result says that the only way a regular solution might blow-up is if the $L^\infty$-norm of $1/\rho$ blows-up. Theorem \ref{Cons} is essentially an adaptation to the whole space of the blow-up criterion proved in Constantin \textit{et al} (see Theorem $1.1.$ from \cite{Constantin}) in the case of the torus. We refer the reader to the Appendix for a sketch of the proof.

The objective of the rest of this section and section \ref{s3.1}  is to show how to obtain a new bound for the $L^\infty$-norm of $1/\rho$ by analysing a new quantity that we call effective pressure. Consider a pair $(\rho_0,u_0)$ verifying the hypothesis stated in Theorem \ref{theo1} and let us also consider the following sequence:
\[
\begin{aligned}
&\rho_0^n-1=j_n*(\rho_0-1)\;\;\mbox{and}\;\;v_0^n=j_n*v_0,
\end{aligned}
\]
with $j_{n}$ a regularizing kernel, $j_{n}(y)=nj(ny)$ with $0\leq j\leq1$,
$\int_{\mathbb{R}}j(y)dy=1$, $j\in C^{\infty}(\mathbb{R})$ and $\mbox{supp}j\subset[-2,2]$. We deduce that
$(\rho_{0}^{n}-1,v_{0}^{n})$ belong to all Sobolev spaces $H^{s}(\mathbb{R})$ with $s\geq 5/2$ and that:
\begin{equation}
0<c\leq\rho_{0}^{n}\leq M<+\infty.
\label{good}
\end{equation}
Also, we consider%
\[
u_{0}^{n}=v_{0}^{n}-\partial_{x}\varphi(\rho_{0}^{n}),
\]
by composition theorem we know that $\varphi(\rho_{0}^{n})-\varphi(1)$ belongs
to $H^{k}(\mathbb{R})$ for any $k\geq0$. Then we obtain that $u_{0}^{n}\in
H^{k}(\mathbb{R})$ for $k\geq3$. Finally we have for $x>y$ and using
(\ref{BV}):
\[
\frac{v_{0}^{n}(x)-v_{0}^{n}(y)}{x-y}=\int_{\mathbb{R}}(\frac{v_{0}%
(x-z)-v_{0}(y-z)}{x-y})j_{n}(z)dz\leq C.
\]
In particular we deduce that for any $x\in\mathbb{R}$, we have:
\begin{equation}
\partial_{x}v_{0}^{n}(x)\leq C.\label{imp1}%
\end{equation}
Using the Theorem (\ref{Cons}), we deduce that there exists a strong solution
$(\rho_{n},u_{n})$ on $(0,T_{n}^{\ast})$ with $n\in\mathbb{N}$. We are going
to prove that $T_{n}^{\ast}=+\infty$ and uniform estimates on $(\rho_{n}%
,u_{n})_{n\in\mathbb{N}}$ on the time interval $\mathbb{R}^{+}$. The goal now
is to apply the blow-up criterion of the Theorem \ref{Cons}. Let us prove that
for any $t\in(0,T_{n}^{\ast})$:
\[
\Vert\frac{1}{\rho_{n}}(t,\cdot)\Vert_{L^{\infty}}\leq C,
\]
for any $n\in\mathbb{N}$. Let us recall that there exists $C>0$ such that for
any $t>0$ we have:
\begin{equation}
\begin{aligned} &\int_{\R}[\rho_n(t,x)|u_n|^2(t,x)+\Pi(\rho_n(t,x))-\Pi(1)] dx+\int^t_0\int_{\R}\mu(\rho_n(s,x))(\partial_x u_n(s,x))^2 ds dx\leq C, \end{aligned}\label{energy}%
\end{equation}
and:
\begin{equation}
\begin{aligned} &\int_{\R}[\rho_n(t,x)|v_n|^2(t,x)+\Pi(\rho_n(t,x))-\Pi(1)] dx+\int^t_0\int_{\R}\frac{\mu(\rho_n)P'(\rho_n)}{\rho_n} |\partial_x\rho_n(s,x)|^2 ds dx\leq C, \end{aligned}\label{BD}%
\end{equation}
This is due to the fact that it exists $C_{1}>0$ such that:
\[
\begin{aligned}
&\|v_0^n\|_{L^2(\R)}\leq C_1,\;\|\rho_0^n-1\|_{L^2(\R)}\leq C_1\;\;\mbox{and}\;\;\|\p_x\rho_0^n\|_{L^2(\R)}\leq C_1.
\end{aligned}
\]
Combining (\ref{energy}) and (\ref{BD}), we deduce that for $C>0$ large enough
we have for any $t\in(0,T_{n}^{\ast})$:
\begin{equation}
\Vert\rho_{n}(t,\cdot)-1\Vert_{L_{2}^{\gamma}(\mathbb{R})}\leq C,\;\Vert
\sqrt{\rho_{n}}\partial_{x}\varphi(\rho_{n})\Vert\leq C.\label{mBD}%
\end{equation}
We refer to \cite{Lio98} for the definition of Orlicz spaces. Since
$\gamma>\alpha+1$, using (\ref{mBD}) and the Lemma 3.7 from (\cite{Jiu}) we
get for $C>0$ large enough and independent on $n$:
\begin{equation}
\Vert\rho_{n}\Vert_{L^{\infty}([0,T_{n}^{\ast}],L^{\infty})}\leq
C.\label{Linf}%
\end{equation}

\subsection{New effective pressure $y_{n}$ and uniform estimates for $\frac{1}{\rho_n}$}\label{s3.1}

We recall now that the effective velocity $v_{n}$ verifies the momentum
equation of the system (\ref{0.1a}), namely:
\[
\begin{aligned}
&\p_t v_n+u_n\p_x v_n+\p_x F(\rho_n)=0,
\end{aligned}
\]
with:
\[
\partial_{x}F(\rho_{n})=\frac{P^{\prime}(\rho_{n})\rho_{n}}{\mu(\rho_{n}%
)}(v_{n}-u_{n}).
\]
Let us set now $w_{n}=\partial_{x}v_{n}$, we observe that $w_{n}$ satisfies
the following equation:
\[
\begin{aligned}
&\p_t w_n+u_n\p_x w_n+\p_x u_n w_n+\frac{P'(\rho_n)\rho_n}{\mu(\rho_n)} w_n-\frac{P'(\rho_n)\rho_n}{\mu(\rho_n)}\p_x u_n+\p_x(\frac{P'(\rho_n)\rho_n}{\mu(\rho_n)})(v_n-u_n)=0.
\end{aligned}
\]
If we set $F_{1}(\rho)=\frac{P^{\prime}(\rho)\rho}{\mu(\rho)}$, we have:
\[
\begin{aligned}
&\p_t w_n+u_n\p_x w_n+\p_x u_n w_n+F_1(\rho_n) w_n-F_1(\rho_n)\p_x u_n+F_1' (\rho_n)\frac{\rho_n^2}{\mu(\rho_n)}(v_n-u_n)^2=0.
\end{aligned}
\]
Let us multiply the previous equation by $\frac{1}{\rho_{n}}$, we get then:
\[
\begin{aligned}
&\p_t (\frac{w_n}{\rho_n})+u_n\p_x (\frac{w_n}{\rho_n})+F_1(\rho_n) \frac{w_n}{\rho_n}-\frac{F_1(\rho_n)}{\rho_n}\p_x u_n+F_1' (\rho_n)\frac{\rho_n}{\mu(\rho_n)}(v_n-u_n)^2=0.
\end{aligned}
\]
We set now $y_{n}=\frac{w_{n}}{\rho_{n}}+F_{2}(\rho_{n})$ with $\rho_{n}%
F_{2}^{\prime}(\rho_{n})=\frac{F_{1}(\rho_{n})}{\rho_{n}}$, we obtain then:
\begin{equation}
\begin{aligned} &\partial_t y_n+u_n\partial_x y_n+F_1(\rho_n) y_n-F_1(\rho_n)F_2(\rho_n)+F_1' (\rho_n)\frac{\rho_n}{\mu(\rho_n)}(v_n-u_n)^2=0. \end{aligned}\label{cru}%
\end{equation}
We recall now that $P(\rho)=\rho_{n}^{\gamma}$, $\mu(\rho_{n})=\rho
_{n}^{\alpha}$ and we get:
\begin{equation}%
\begin{cases}
\begin{aligned} &F_2(\rho_n)=\frac{\gamma}{\gamma-\alpha-1}\rho_n^{\gamma-\alpha-1}\;\;\;\mbox{if}\;\gamma-\alpha-1\ne 0\\ &F_2(\rho_n)=\gamma\ln\rho_n\;\;\;\mbox{if}\;\gamma=\alpha+1\\ &F_1(\rho_n)=\gamma\rho_n^{\gamma-\alpha}. \end{aligned}
\end{cases}
\label{acalcul}%
\end{equation}
Now since $y_{n}$ is continuous on $\mathbb{R}^{+}\times\mathbb{R}$ (indeed we
recall that the solution $(\rho_{n},u_{n})$ is regular) and $\lim
_{x\rightarrow\pm\infty}y_{n}(t,x)=F_{2}(1)$, we deduce that $y_{n}(t,\cdot)$
has a maximum for every $t\geq0$ and thus the function $y_{n}^{M}$ defined as
follows:
\[
y_{M}^{n}(t)=\max_{x\in\mathbb{R}}y_{n}(t,x).
\]
makes sense. Furthermore $y_{n}^{M}$ is Lipschitz continuous on any interval
$[0,T]$ with $T\in(0,T_{n}^{\ast})$. Indeed from the triangular inequality we
have for $(t_{1},t_{2})\in(0,T_{n}^{\ast})$:
\[
|y_{M}^{n}(t_{1})-y_{M}^{n}(t_{1})|\leq\max_{x\in\mathbb{R}}|y_{n}%
(t_{1},x)-y_{,}(t_{2},x)|\leq\Vert\partial_{s}y_{n}\Vert_{L^{\infty}%
([t_{1},t_{2}],L^{\infty})}|t_{1}-t_{2}|.
\]
According to Rademacher theorem, $y_{n}^{M}$ is differentiable almost
everywhere on $[0,T_{n}^{\ast})$. Furthermore there exists for each
$t\in\lbrack0,T_{n}^{\ast})$ a point $x_{t}^{n}$ such that:
\[
y_{M}^{n}(t)=y(t,x_{t}^{n}).
\]
We are going to verify now that for almost all $t\in(0,T_{n}^{\ast})$ we have
$(y_{M}^{n})^{\prime}(t)=\partial_{t}y_{n}(t,x_{t}^{n})$. Indeed we have:
\[
\begin{aligned}
(y^n_M)'(t)&=\lim_{h\rightarrow 0^+}\frac{y^n_M(t+h)-y^n_M(t)}{h}=\lim_{h\rightarrow 0^+}\frac{y_n(t+h,x_{t+h}^n)-y^n(t,x_t^n)}{h}\\
&\geq \lim_{h\rightarrow 0^+}\frac{y_n(t+h,x_{t}^n)-y^n(t,x_t^n)}{h}=\p_t y_n(t,x_t^n).
\end{aligned}
\]
Similarly, we have:
\[
\begin{aligned}
(y^n_M)'(t)&=\lim_{h\rightarrow 0^+}\frac{y^n_M(t)-y^n_M(t-h)}{h}=\lim_{h\rightarrow 0^+}\frac{y_n(t,x_{t}^n)-y^n(t-h,x_{t-h}^n)}{h}\\
&\leq \lim_{h\rightarrow 0^+}\frac{y_n(t,x_{t}^n)-y^n(t-h,x_t^n)}{h}=\p_t y_n(t,x_t^n).
\end{aligned}
\]
We deduce from (\ref{cru}) using the fact that $\partial_{x}y_{n}(t,x_{t}%
^{n})=0$ since $y_{n}(t,\cdot)$ reaches its maximum in $x_{t}^{n}$ and that
for all $t\in(0,T_{n}^{\ast})$ we have:
\begin{equation}
\begin{aligned} &\partial_t y_M^n(t)+F_1(\rho_n)(t,x_t^n)y_M^n(t)=F_1(\rho_n)F_2(\rho_n)(t,x_t^n)-F_1' (\rho_n)\frac{\rho_n}{\mu(\rho_n)}(v_n-u_n)^2(t,x_t^n). \end{aligned}\label{cru1}%
\end{equation}
Basic computations give now:
\begin{equation}%
\begin{cases}
\begin{aligned} &F_1(\rho)F_2(\rho)=\frac{\gamma^2}{\gamma-\alpha-1}\rho^{2\gamma-2\alpha-1}\;\;\;\mbox{if}\;\;\;\gamma\ne\alpha+1\\ &F_1(\rho)F_2(\rho)=\gamma^2 \ln\rho \rho^{\gamma-\alpha}\;\;\;\mbox{if}\;\;\;\gamma=\alpha+1\\ &\frac{F'_1(\rho)\rho}{\mu(\rho)}=\gamma(\gamma-\alpha)\rho^{\gamma-2\alpha}\\ &F_1(\rho)=\gamma\rho^{\gamma-\alpha}. \end{aligned}
\end{cases}
\label{calcul}%
\end{equation}
We recall that we have $\gamma\geq \alpha$ such that using (\ref{cru1}) and
(\ref{calcul}) we get that for $\gamma\ne\alpha+1$:
\begin{equation}
\begin{aligned} &\partial_t y_M^n(t)+F_1(\rho_n)(t,x_t^n)y_M^n(t)\leq \max(0,\frac{\gamma^2}{\gamma-\alpha-1})\|\rho_n(t,\cdot)\|_{L^\infty}^{2\gamma-2\alpha-1}. \end{aligned}\label{cru2}%
\end{equation}
From (\ref{calcul}) and (\ref{cru2}), we get that for any $t\in\lbrack0,T_{n}^{\ast})$ one
has with $C_\gamma=\max(0,\frac{\gamma^2}{\gamma-\alpha-1})$:
\begin{equation}
\partial_{t}(y_{M}^{n}(t)e^{\gamma\int_{0}^{t}\rho_{n}^{\gamma-\alpha
}(s,x_{s}^{n})ds})\leq C_\gamma\Vert\rho
_{n}(t,\cdot)\Vert_{L^{\infty}}^{2\gamma-2\alpha-1}e^{\gamma\int_{0}^{t}%
\rho_{n}^{\gamma-\alpha}(s,x_{s}^{n})ds}.
\end{equation}
It yields for any $t\in(0,T_{n}^{\ast})$:
\begin{equation}
\begin{aligned} &y_M^n(t)\leq e^{-\gamma\int^t_0\rho_n^{\gamma-\alpha}(s,x_s^n) ds}y_M^n(0)\\ &\hspace{3cm}+C_\gamma\int^t_0 \|\rho_n(t,\cdot)\|_{L^\infty}^{2\gamma-2\alpha-1} e^{-\gamma\int^t_s\rho_n^{\gamma-\alpha}(s',x_s^n) ds'} ds . \end{aligned}
\label{accru2}
\end{equation}
From (\ref{imp1}) and (\ref{good}), we deduce that for any $n\in\mathbb{N}$ we have for any $x\in\R$ and $\gamma\ne\alpha+1$:
\begin{equation}
\begin{aligned}
y^n(0,x)&\leq \frac{\max(0,C)}{c}+\frac{\gamma}{\gamma-\alpha-1} M^{\gamma-\alpha-1}=C_1.
\end{aligned}
\label{gini}
\end{equation}
We obtain now from (\ref{accru2}), (\ref{gini}) and since $\rho_n$ is positive:
\begin{equation}
\begin{aligned} &y_M^n(t)\leq C_1+C_\gamma \int^t_0 \|\rho_n(t,\cdot)\|_{L^\infty}^{2\gamma-2\alpha-1} ds . \end{aligned}\label{cru4}%
\end{equation}
Combining (\ref{Linf}) and (\ref{cru4}), we deduce that for any $t\in
(0,T_{n}^{\ast})$ we have: 
\begin{equation}
y_{M}^{n}(t)\leq C(t),\label{pointfin10}%
\end{equation}
with $C$ a continuous function on $\mathbb{R}^{+}$ when $\gamma\ne\alpha+1$. 
From (\ref{pointfin10}), it yields for any $t\in(0,T_{n}^{\ast})$ and
$x\in\mathbb{R}$ when $\gamma\ne\alpha+1$:
\begin{equation}
\frac{\partial_{x}v_{n}(t,x)}{\rho_n(t,x)}\leq C(t)+\frac{\gamma}{\alpha+1-\gamma}\rho_n^{\gamma-\alpha-1}(t,x),\label{vvpointfin1}%
\end{equation}
with $C$ a continuous function on $\mathbb{R}^{+}$. Next we recall that we
have:
\[
\partial_{t}(\frac{1}{\rho_{n}})+u_{n}\partial_{x}(\frac{1}{\rho_{n}}%
)-\frac{1}{\rho_{n}}\partial_{x}u_{n}=0.
\]
We can rewrite the equation as follows:
\[
\begin{aligned}
&\p_t(\frac{1}{\rho_n})+u_n\p_x(\frac{1}{\rho_n})-\frac{1}{\rho_n}\p_x v_n-\frac{\mu(\rho_n)}{\rho_n}\p_{xx}(\frac{1}{\rho_n})-\frac{1}{\rho_n}\p_x\mu(\rho_n)\p_x(\frac{1}{\rho_n})=0.
\end{aligned}
\]
Using again a maximum principle and following the same arguments as 
previously, we set now:
\[
z_{n}(t)=\sup_{x\in\mathbb{R}}\frac{1}{\rho_{n}}(t,x)=\frac{1}{\rho_{n}%
}(t,x_{t}^{n}).
\]
We have then:
\[
\begin{aligned}
&\p_t z_n(t)=\frac{\mu(\rho_n)}{\rho_n}\p_{xx}(\frac{1}{\rho_n})(t,x^n_t)+\frac{1}{\rho_n}\p_x v_n(t,x^n_t).
\end{aligned}
\]
From (\ref{vvpointfin1}) and since $\partial_{xx}(\frac{1}{\rho_{n}}%
)(t,x_{t}^{n})\leq0$ (indeed $x_{t}^{n}$ is a point where $\frac{1}{\rho_{n}}$
reaches its maximum) we deduce that:
\begin{equation}
\begin{aligned}
\p_t z_n(t)&\leq C(t)+\frac{\gamma}{\alpha+1-\gamma}\rho_n^{\gamma-\alpha-1}(t,x^n_t)\\
&\leq  C(t)+\frac{\gamma}{\alpha+1-\gamma} z_n(t)^{\alpha+1-\gamma}.
\end{aligned}
\end{equation}
Using Gronwall lemma, it implies that there exists a continuous function $C_{2}$
on $\mathbb{R}^{+}$ such that for any $t\in(0,T_{n}^{\ast})$ we have:
\[
z_{n}(t)\leq C_{2}(t).
\]
This implies that for any $t\in(0,T_{n}^{\ast})$ we get:
\begin{equation}
\Vert\frac{1}{\rho_{n}}(t,\cdot)\Vert_{L^{\infty}}\leq C_{2}(t).\label{vide}%
\end{equation}
Combining the blow-up criterion in Theorem \ref{Cons} and (\ref{vide}), we
obtain that $T_{n}^{\ast}=+\infty$ and for any $t>0$:
\begin{equation}
\Vert\frac{1}{\rho_{n}}(t,\cdot)\Vert_{L^{\infty}}\leq C_{2}%
(t),\label{loin_de_vide}%
\end{equation}
with $C_{2}$ a continuous function on $\mathbb{R}^{+}$.
From (\ref{vvpointfin1}), (\ref{Linf})
and (\ref{loin_de_vide}), we get again for any $t\in(0,T_{n}^{\ast})$ and
$x\in\mathbb{R}$ when $\gamma\ne\alpha+1$:
\begin{equation}
\partial_{x}v_{n}(t,x)\leq C_{1}(t),\label{pointfin1}%
\end{equation}
with $C_1$ a continuous increasing function. We can easily prove similar results for $\gamma=\alpha+1$.
\subsection{Estimates \`{a} la Hoff}

In the sequel for simplifying the notation we drop the index $n$. Introducing
the convective derivative
\[
\dot{u}=\partial_{t}u+u\partial_{x}u,
\]
we rewrite the momentum equation as%
\[
\rho\dot{u}-\partial_{x}\left(  \rho^{\alpha}u_{x}\right)  +\partial_{x}%
\rho^{\gamma}=0.
\]
Let us observe that:
\begin{align}
-\int_{\mathbb{R}}\partial_{x}\left(  \rho^{\alpha}\partial_{x}u\right)
\partial_{t}u  &  =\int_{\mathbb{R}}\rho^{\alpha}\partial_{x}u\partial
_{xt}^{2}u=\frac{1}{2}\int_{\mathbb{R}}\rho^{\alpha}\partial_{t}\left(
(\partial_{x}u)^{2}\right) \nonumber\\
&  =\frac{1}{2}\frac{d}{dt}\int_{\mathbb{R}}\rho^{\alpha}\left(  \partial
_{x}u\right)  ^{2}-\frac{1}{2}\int_{\mathbb{R}}\partial_{t}\rho^{\alpha
}(\partial_{x}u)^{2}. \label{estim1}%
\end{align}
Next, we see that:%
\begin{align*}
-\int_{\mathbb{R}}\partial_{x}\left(  \rho^{\alpha}\partial_{x}u\right)
u\partial_{x}u  &  =-\int_{\mathbb{R}}u\partial_{x}\rho^{\alpha}(\partial
_{x}u)^{2}-\int_{\mathbb{R}}\rho^{\alpha}u\partial_{xx}^{2}u\partial_{x}u\\
&  =-\int_{\mathbb{R}}u\partial_{x}\rho^{\alpha}(\partial_{x}u)^{2}+\frac
{1}{2}\int_{\mathbb{R}}\partial_{x}\left(  u\rho^{\alpha}\right)
(\partial_{x}u)^{2}\\
&  =-\int_{\mathbb{R}}u\partial_{x}\rho^{\alpha}(\partial_{x}u)^{2}+\frac
{1}{2}\int_{\mathbb{R}}\rho^{\alpha}(\partial_{x}u)^{3}+\frac{1}{2}%
\int_{\mathbb{R}}u\partial_{x}\rho^{\alpha}(\partial_{x}u)^{2}\\
&  =-\frac{1}{2}\int_{\mathbb{R}}u\partial_{x}\rho^{\alpha}(\partial_{x}%
u)^{2}+\frac{1}{2}\int_{\mathbb{R}}\rho^{\alpha}(\partial_{x}u)^{3}.
\end{align*}
Thus, we gather that:%
\begin{align*}
-\int_{\mathbb{R}}\partial_{x}\left(  \rho^{\alpha}\partial_{x}u\right)
\dot{u}  &  =\frac{1}{2}\frac{d}{dt}\int_{\mathbb{R}}\rho^{\alpha}\left(
\partial_{x}u\right)  ^{2}-\frac{1}{2}\int_{\mathbb{R}}\partial_{t}%
\rho^{\alpha}(\partial_{x}u)^{2}-\frac{1}{2}\int_{\mathbb{R}}u\partial_{x}%
\rho^{\alpha}(\partial_{x}u)^{2}+\frac{1}{2}\int_{\mathbb{R}}\rho^{\alpha
}(\partial_{x}u)^{3}\\
&  =\frac{1}{2}\frac{d}{dt}\int_{\mathbb{R}}\rho^{\alpha}\left(  \partial
_{x}u\right)  ^{2}+\frac{1+\alpha}{2}\int_{\mathbb{R}}\rho^{\alpha}%
(\partial_{x}u)^{3}.
\end{align*}
Moreover, we see that:%
\begin{align*}
\int_{\mathbb{R}}\partial_{x}\rho^{\gamma}\left(  \partial_{t}u+u\partial
_{x}u\right)   &  =-\int_{\mathbb{R}}\rho^{\gamma}\partial_{tx}u+\int
_{\mathbb{R}}u\partial_{x}\rho^{\gamma}\partial_{x}u\\
&  =-\frac{d}{dt}\int_{\mathbb{R}}\rho^{\gamma}\partial_{x}u+\int_{\mathbb{R}%
}\partial_{t}\rho^{\gamma}\partial_{x}u+\int_{\mathbb{R}}u\partial_{x}%
\rho^{\gamma}\partial_{x}u\\
&  =-\frac{d}{dt}\int_{\mathbb{R}}\rho^{\gamma}\partial_{x}u-\gamma
\int_{\mathbb{R}}\rho^{\gamma}(\partial_{x}u)^{2}.
\end{align*}
Multiplying the momentum equation with $\dot{u}$ yields:%
\begin{equation}
\int_{\mathbb{R}}\rho\dot{u}^{2}+\frac{d}{dt}\left\{  \frac{1}{2}%
\int_{\mathbb{R}}\rho^{\alpha}\left(  \partial_{x}u\right)  ^{2}%
-\int_{\mathbb{R}}\rho^{\gamma}\partial_{x}u\right\}  =-\frac{1+\alpha}{2}%
\int_{\mathbb{R}}\rho^{\alpha}(\partial_{x}u)^{3}+\gamma\int_{\mathbb{R}}%
\rho^{\gamma}(\partial_{x}u)^{2}.\label{impo}%
\end{equation}
Let us multiply the previous estimate by $\sigma\left(  t\right)  =\min( 1,t )
$ and integrate in time on $[0,t]$ with $t>0$, we have then:%
\begin{align*}
&  \frac{\sigma\left(  t\right)  }{2}\int_{\mathbb{R}}\rho^{\alpha}\left(
t\right)  \left(  \partial_{x}u\right)  ^{2}\left(  t\right)  +\int_{0}%
^{t}\int_{\mathbb{R}}\sigma\rho\dot{u}^{2}\\
&  =\sigma\left(  t\right)  \int_{\mathbb{R}}\left(  \rho^{\gamma}-1\right)
\partial_{x}u+\int_{0}^{\min\left\{  1,t\right\}  }\int_{\mathbb{R}}\left[
\frac{1}{2}\rho^{\alpha}\left(  \partial_{x}u\right)  ^{2}-\left(
\rho^{\gamma}-1\right)  \partial_{x}u\right] \\
&  -\frac{1+\alpha}{2}\int_{0}^{t}\int_{\mathbb{R}}\sigma\rho^{\alpha
}(\partial_{x}u)^{3}+\gamma\int_{0}^{t}\int_{\mathbb{R}}\sigma\rho^{\gamma
}(\partial_{x}u)^{2}.
\end{align*}
Let us denote by:
\[
A\left(  \rho,u\right)  \left(  t\right)  =\frac{\sigma\left(  t\right)  }%
{2}\int_{\mathbb{R}}\rho^{\alpha}\left(  t\right)  \left(  \partial
_{x}u\right)  ^{2}\left(  t\right)  +\int_{0}^{t}\int_{\mathbb{R}}\sigma
\rho\dot{u}^{2}.
\]
Let us observe that using (\ref{energy}), (\ref{Linf}) and (\ref{loin_de_vide}%
) we have:%
\begin{equation}
\begin{aligned} \sigma\left( t\right) \int_{\mathbb{R}}\left( \rho^{\gamma}-1\right) \partial_{x}u & \leq \sqrt{\sigma(t)} \left\Vert \frac{\rho^{\gamma}-1}{\rho^{\frac{\alpha}{2}}}\right\Vert _{L_{t}^{\infty}L^{2}}\left( \int_{\mathbb{R}}\sigma\left( t\right) \rho^{\alpha}\left( t\right) (\partial_{x}u)^{2}\left( t\right) \right) ^{\frac{1}{2}}\\ & \leq C\left( t\right) \left\Vert \frac{\rho^{\gamma}-1}{\rho^{\frac {\alpha}{2}}}\right\Vert _{L_{t}^{\infty}L^{2}}^{2}+\frac{1}{4}\int _{\mathbb{R}}\sigma\left( t\right) \rho^{\alpha}\left( t\right) (\partial_{x}u)^{2}\left( t\right) \\ & \leq C_1\left( t\right) +\frac{1}{4}\int_{\mathbb{R}}\sigma\left( t\right) \rho^{\alpha}\left( t\right) (\partial_{x}u)^{2}\left( t\right), \end{aligned}\label{impo1}%
\end{equation}
with $C$ and $C_{1}$ continuous on $\mathbb{R}^{+}$. Next, we see that owing
to the estimate (\ref{energy}), (\ref{Linf}) and (\ref{loin_de_vide}), we have
that:%
\begin{equation}
\int_{0}^{\min\left\{  1,t\right\}  }\int_{\mathbb{R}}\left[  \frac{1}{2}%
\rho^{\alpha}\left(  \partial_{x}u\right)  ^{2}-\left(  \rho^{\gamma
}-1\right)  \partial_{x}u\right]  +\gamma\int_{0}^{t}\int_{\mathbb{R}}%
\sigma\rho^{\gamma}(\partial_{x}u)^{2}\leq C_{2}\left(  t\right)
,\label{impo2}%
\end{equation}
with $C_{2}$ a continuous function on $\mathbb{R}^{+}$. Combining
(\ref{impo}), (\ref{impo1}) and (\ref{impo2}) , we thus get for all $t\geq0$:%
\begin{align*}
A\left(  \rho,u\right)  \left(  t\right)   &  \leq C\left(  t\right)
+\frac{1}{4}\int_{\mathbb{R}}\sigma\left(  t\right)  \rho^{\alpha}\left(
t\right)  (\partial_{x}u)^{2}\left(  t\right)  -\frac{1+\alpha}{2}\int_{0}%
^{t}\int_{\mathbb{R}}\sigma\rho^{\alpha}(\partial_{x}u)^{3}\\
&  \leq C_{3}\left(  t\right)  +\frac{1}{2}A\left(  \rho,u\right)  \left(
t\right)  -\frac{1+\alpha}{2}\int_{0}^{t}\int_{\mathbb{R}}\sigma\rho^{\alpha
}(\partial_{x}u)^{3}%
\end{align*}
with $C_{3}$ a continuous fonction on $\mathbb{R}^{+}$. Consequently it
yields:%
\[
A\left(  \rho,u\right)  \left(  t\right)  \leq C\left(  t\right)
+(1+\alpha)\int_{0}^{t}\int_{\mathbb{R}}\sigma\rho^{\alpha}(\partial_{x}u)^{3}%
\]
which also implies that ($C$ can be choseen to be increasing in $t$):
\begin{equation}
\sup_{\tau\in\left[  0,t\right]  }A\left(  \rho,u\right)  \left(  \tau\right)
\leq C\left(  t\right)  +(1+\alpha)\int_{0}^{t}\int_{\mathbb{R}}\sigma
\rho^{\alpha}(\partial_{x}u)^{3} \label{(relatie_A)}%
\end{equation}
Let us observe that for all $\varepsilon>0$ we have using Gagliardo-Nirenberg
inequality (\ref{energy}) and (\ref{Linf}):%
\begin{align}
&  \int_{0}^{t}\sigma^{\frac{1}{2}}\left(  \tau\right)  \left\Vert \left(
\rho^{\alpha}\partial_{x}u-\rho^{\gamma}\right)  \left(  \tau\right)
\right\Vert _{L^{\infty}}^{2}\nonumber  \leq2\int_{0}^{t}\sigma^{\frac{1}{2}}\left(  \tau\right)  \left\Vert
\left(  \rho^{\alpha}\partial_{x}u-(\rho^{\gamma}-1)\right)  \left(
\tau\right)  \right\Vert _{L^{\infty}}^{2}+2t\\
&  \leq2\int_{0}^{t}\sigma^{\frac{1}{2}}\left(  \tau\right)  \left\Vert
\left(  \rho^{\alpha}\partial_{x}u-(\rho^{\gamma}-1)\right)  \left(
\tau\right)  \right\Vert _{L^{2}}\left\Vert \partial_{x}\left(  \rho^{\alpha
}\partial_{x}u-\rho^{\gamma}\right)  \left(  \tau\right)  \right\Vert _{L^{2}%
}+2t\nonumber\\
&  \leq C_{\varepsilon}\int_{0}^{t}\left\Vert \left(  \rho^{\alpha}%
\partial_{x}u-(\rho^{\gamma}-1)\right)  \left(  \tau\right)  \right\Vert
_{L^{2}}^{2}+\varepsilon\int_{0}^{t}\sigma\left(  \tau\right)  \left\Vert
\partial_{x}\left(  \rho^{\alpha}\partial_{x}u-\rho^{\gamma}\right)  \left(
\tau\right)  \right\Vert _{L^{2}}^{2}+2t\nonumber\\
& \leq C_{\varepsilon}\int_{0}^{t}\left\Vert \left(  \rho^{\alpha}\partial
_{x}u-(\rho^{\gamma}-1)\right)  \left(  \tau\right)  \right\Vert _{L^{2}}%
^{2}+\varepsilon\int_{0}^{t}\sigma\left(  \tau\right)  \left\Vert \rho\dot
{u}\left(  \tau\right)  \right\Vert _{L^{2}}^{2}+2t\nonumber\\
&  \leq C\left(  t,\varepsilon\right)  +\varepsilon\|\rho\|_{L^{\infty
}([0,t],L^{\infty})}^{\frac{1}{2}}A\left(  \rho,u\right)  \left(  t\right) \\
&  \leq C\left(  t,\varepsilon\right)  +\varepsilon C_{0} A\left(
\rho,u\right)  \left(  t\right)  , \label{Linfty}%
\end{align}
with $C$ a continuous function on $\mathbb{R}^{+}$. We are going now to
estimate the last term of $\left(  \text{\ref{(relatie_A)}}\right)  $ and
using (\ref{energy}), (\ref{Linf}), (\ref{loin_de_vide}) and $\left(
\text{\ref{Linfty}}\right)  $ with $\varepsilon=1/(2\left(  1+\alpha\right)
C_{0})$ we obtain that:%
\begin{align}
&  \int_{0}^{t}\int_{\mathbb{R}}\sigma\rho^{\alpha}(\partial_{x}%
u)^{3}\nonumber =\int_{0}^{t}\int_{\mathbb{R}}\sigma(\partial_{x}u)^{2}(\rho^{\alpha
}\partial_{x}u-\rho^{\gamma})+\int_{0}^{t}\int_{\mathbb{R}}\sigma\rho^{\gamma
}(\partial_{x}u)^{2}\\
&  \leq\int_{0}^{t}\left(  \sigma^{\frac{1}{4}}\left\Vert \left(  \rho
^{\alpha}\partial_{x}u-\rho^{\gamma}\right)  \left(  \tau\right)  \right\Vert
_{L^{\infty}}\sigma^{\frac{3}{4}}\int_{\mathbb{R}}(\partial_{x}u)^{2}\left(
\tau\right)  \right)  d\tau+\int_{0}^{t}\int_{\mathbb{R}}\sigma\rho^{\gamma
}(\partial_{x}u)^{2}\nonumber\\
&  \leq C\left(  t\right)  +\int_{0}^{t}\sigma^{\frac{1}{2}}\left(
\tau\right)  \left\Vert \left(  \rho^{\alpha}\partial_{x}u-\rho^{\gamma
}\right)  \left(  \tau\right)  \right\Vert _{L^{\infty}}^{2}+\int_{0}%
^{t}\sigma^{\frac{3}{2}}\left(  \tau\right)  \left(  \int_{\mathbb{R}%
}(\partial_{x}u)^{2}\left(  \tau\right)  dx\right)  ^{2}\nonumber\\
&  \leq C\left(  t\right)  +\frac{1}{2\left(  1+\alpha\right)  }A\left(
\rho,u\right)  \left(  t\right)  +\int_{0}^{t}\left\Vert \frac{1}{\rho\left(
\tau\right)  }\right\Vert _{L^{\infty}}^{2\alpha}\sigma^{\frac{3}{2}}\left(
\tau\right)  (\int_{\mathbb{R}}\rho^{\alpha}(\partial_{x}u)^{2}\left(
\tau\right)  dx )^{2}\nonumber\\
&  \leq C\left(  t\right)  +\frac{1}{2\left(  1+\alpha\right)  }A\left(
\rho,u\right)  \left(  t\right)  +C_{1}\left(  t\right)  \int_{0}^{t}%
\sigma\left(  \tau\right)  \int_{\mathbb{R}}\rho^{\alpha}(\partial_{x}%
u)^{2}\left(  \tau\right)  dx \int_{\mathbb{R}}(\rho^{\alpha}\partial
_{x}u)^{2}\left(  \tau\right)  d x\nonumber\\
&  \leq C\left(  t\right)  +\frac{1}{2\left(  1+\alpha\right)  }A\left(
\rho,u\right)  \left(  t\right)  +2C_{1}\left(  t\right)  \int_{0}^{t}A\left(
\rho,u\right)  \left(  \tau\right)  \int_{\mathbb{R}}(\rho^{\alpha}%
\partial_{x}u)^{2}\left(  \tau\right)  d\tau,\label{relatie_A2}%
\end{align}
with $C$ and $C_{1}$ continuous increasing functions. Finally, putting
together $\left(  \text{\ref{(relatie_A)}}\right)  $ and $\left(
\text{\ref{relatie_A2}}\right)  $ we get that%
\[
\sup_{\tau\in\left[  0,t\right]  }A\left(  \rho,u\right)  \left(  \tau\right)
\leq C_{2}\left(  t\right)  +C_{2}\left(  t\right)  \int_{0}^{t}A\left(
\rho,u\right)  \left(  \tau\right)  \int_{\mathbb{R}}(\rho^{\alpha}%
\partial_{x}u)^{2}\left(  \tau\right)  d\tau,
\]
with $C_{2}$ an increasing continuous function. Using Gronwall's lemma and
(\ref{energy}) leads to%
\begin{equation}
\sup_{\tau\in\left[  0,t\right]  }A\left(  \rho,u\right)  \left(  \tau\right)
\leq C\left(  t\right)  , \label{estimate_for_A(rho,u)}%
\end{equation}
with $C$ an increasing continuous function. The control over $A\left(
\rho,u\right)  $ and $\left(  \text{\ref{Linfty}}\right)  $ yields%
\[
\int_{0}^{t}\sigma^{\frac{1}{2}}\left(  \tau\right)  \left\Vert \left(
\rho^{\alpha}\partial_{x}u-\rho^{\gamma}\right)  \left(  \tau\right)
\right\Vert _{L^{\infty}}^{2}d\tau\leq C\left(  t\right) ,
\]
and consequently we get using in addition (\ref{Linf}):%
\begin{equation}
\int_{0}^{t}\sigma^{\frac{1}{2}}\left(  \tau\right)  \left\Vert \partial
_{x}u\left(  \tau\right)  \right\Vert _{L^{\infty}}^{2}d\tau\leq C\left(
t\right) .\label{Lip}%
\end{equation}
The last inequality also provides an estimate in $L_{t}^{1}\left(  L^{\infty
}\right)  $ of $\partial_{x}u$ for any $t>0$ using Cauchy-Schwarz inequality:%
\[
\int_{0}^{t}\left\Vert \partial_{x}u\left(  \tau\right)  \right\Vert
_{L^{\infty}}d\tau\leq\left(  \int_{0}^{t}\sigma^{-\frac{1}{2}}\left(
\tau\right)  d\tau\right)  ^{\frac{1}{2}}\left(  \int_{0}^{t}\sigma^{\frac
{1}{2}}\left(  \tau\right)  \left\Vert \partial_{x}u\left(  \tau\right)
\right\Vert _{L^{\infty}}^{2}\right)  ^{\frac{1}{2}}\leq C\left(  t\right) .
\]
Next, we aim at obtaining estimate for the $L^{2}$-norm of $\partial_{x}%
\dot{u}$. This will be useful in order to recover regularity properties of $u
$. The idea is to apply the operator $\partial_{t}+u\partial_{x}$ to the
velocity's equation:%
\[
\left(  \partial_{t}+u\partial_{x}\right)  (\rho\dot{u})-(\partial
_{t}+u\partial_{x})\partial_{x}\left(  \rho^{\alpha}u_{x}\right)  +\left(
\partial_{t}P(\rho)+u\partial_{x}P(\rho)\right)  =0
\]
and to test it with $\min\{1,t\}\dot{u}$. We begin by observing that
\[
\int_{\mathbb{R}}\left(  \rho\dot{u}\right)  _{t}\dot{u}=\int_{\mathbb{R}}%
\rho_{t}\dot{u}^{2}+\frac{1}{2}\int_{\mathbb{R}}\rho\frac{d\dot{u}^{2}}%
{dt}=\frac{1}{2}\frac{d}{dt}\int_{\mathbb{R}}\rho\dot{u}^{2}+\frac{1}{2}%
\int_{\mathbb{R}}\rho_{t}\dot{u}^{2}.
\]
We remark that:%
\[
\int_{\mathbb{R}}u\partial_{x}\left(  \rho\dot{u}\right)  \dot{u}%
=-\int_{\mathbb{R}}\rho\dot{u}\partial_{x}\left(  u\dot{u}\right)
=-\int_{\mathbb{R}}\partial_{x}u\rho\dot{u}^{2}+\frac{1}{2}\int_{\mathbb{R}%
}\left(  \rho u\right)  _{x}\dot{u}^{2}.
\]
Summing the above two relations yields:%
\begin{equation}
\int_{\mathbb{R}}\left(  \partial_{t}+u\partial_{x}\right)  (\rho\dot{u}%
)\dot{u}=\frac{1}{2}\frac{d}{dt}\int_{\mathbb{R}}\rho\dot{u}^{2}%
-\int_{\mathbb{R}}\partial_{x}u\rho\dot{u}^{2}. \label{Higer_estimates_1}%
\end{equation}
Next, we take a look at the second term:
\begin{align}
\label{Higer_estimates_2a}-\int_{\mathbb{R}}(\partial_{t}+u\partial
_{x})\partial_{x}\left(  \rho^{\alpha}\partial_{x}u\right)  \dot{u}  &
=\int_{\mathbb{R}}\partial_{t}\rho^{\alpha}\partial_{x}u\partial_{x}\dot
{u}+\int_{\mathbb{R}}\rho^{\alpha}\partial_{x}u_{t}\partial_{x}\dot{u}%
+\int_{\mathbb{R}}\partial_{x}(\rho^{\alpha}\partial_{x}u)\partial_{x}%
(u\dot{u})\nonumber\\
\end{align}
Let us treat separately the last term appearing in the above inequality :%
\begin{align}
&  \int_{\mathbb{R}}\partial_{x}(\rho^{\alpha}\partial_{x}u)\partial_{x}%
(u\dot{u})\nonumber\\
&  =\int_{\mathbb{R}}\partial_{x}\rho^{\alpha}(\partial_{x}u)^{2}\dot{u}%
+\int_{\mathbb{R}}u\partial_{x}\rho^{\alpha}\partial_{x}u\partial_{x}\dot
{u}+\int_{\mathbb{R}}\rho^{\alpha}\partial_{xx}^{2}u\partial_{x}u\dot{u}%
+\int_{\mathbb{R}}\rho^{\alpha}u\partial_{xx}^{2}u\partial_{x}\dot
{u}\nonumber\\
&  =\int_{\mathbb{R}}\partial_{x}\rho^{\alpha}(\partial_{x}u)^{2}\dot{u}%
+\int_{\mathbb{R}}u\partial_{x}\rho^{\alpha}\partial_{x}u\partial_{x}\dot
{u}-\frac{1}{2}\int_{\mathbb{R}}(\partial_{x}u)^{2}\partial_{x}(\rho^{\alpha
}\dot{u})+\int_{\mathbb{R}}\rho^{\alpha}\partial_{x}(u\partial_{x}%
u)\partial_{x}\dot{u}-\int_{\mathbb{R}}\left(  \partial_{x}u\right)  ^{2}%
\rho^{\alpha}\partial_{x}\dot{u}\nonumber\\
&  =\frac{1}{2}\int_{\mathbb{R}}\partial_{x}\rho^{\alpha}(\partial_{x}%
u)^{2}\dot{u}+\int_{\mathbb{R}}u\partial_{x}\rho^{\alpha}\partial_{x}%
u\partial_{x}\dot{u}-\frac{3}{2}\int_{\mathbb{R}}(\partial_{x}u)^{2}%
\rho^{\alpha}\partial_{x}\dot{u}+\int_{\mathbb{R}}\rho^{\alpha}\partial
_{x}(u\partial_{x}u)\partial_{x}\dot{u} \label{Higer_estimates_2b}%
\end{align}
Combining the two identities $\left(  \text{\ref{Higer_estimates_2a}}\right)
$ and $\left(  \text{\ref{Higer_estimates_2b}}\right)  $ we get that%
\begin{align}
& -\int_{\mathbb{R}}(\partial_{t}+u\partial_{x})\partial_{x}\left(
\rho^{\alpha}\partial_{x}u\right)  \dot{u} =\int_{\mathbb{R}}\partial_{t}%
\rho^{\alpha}\partial_{x}u\partial_{x}\dot{u}+\int_{\mathbb{R}}u\partial
_{x}\rho^{\alpha}\partial_{x}u\partial_{x}\dot{u}\nonumber\\
&  +\int_{\mathbb{R}}\rho^{\alpha}\partial_{x}u_{t}\partial_{x}\dot{u}%
+\int_{\mathbb{R}}\rho^{\alpha}\partial_{x}(u\partial_{x}u)\partial_{x}\dot{u}
-\frac{3}{2}\int_{\mathbb{R}}(\partial_{x}u)^{2}\rho^{\alpha}\partial_{x}%
\dot{u}+\frac{1}{2}\int_{\mathbb{R}}\partial_{x}\rho^{\alpha}(\partial
_{x}u)^{2}\dot{u}\nonumber\\[2mm]
&  =-\alpha\int_{\mathbb{R}}\rho^{\alpha}(\partial_{x}u)^{2}\partial_{x}%
\dot{u}+\int_{\mathbb{R}}\rho^{\alpha}(\partial_{x}\dot{u})^{2}-\frac{3}%
{2}\int_{\mathbb{R}}(\partial_{x}u)^{2}\rho^{\alpha}\partial_{x}\dot{u}%
+\frac{1}{2}\int_{\mathbb{R}}\partial_{x}\rho^{\alpha}(\partial_{x}u)^{2}%
\dot{u}\nonumber\\
&  =\int_{\mathbb{R}}\rho^{\alpha}(\partial_{x}\dot{u})^{2}-\left(
\alpha+\frac{3}{2}\right)  \int_{\mathbb{R}}\rho^{\alpha}(\partial_{x}%
u)^{2}\partial_{x}\dot{u}+\frac{1}{2}\int_{\mathbb{R}}\partial_{x}\rho
^{\alpha}(\partial_{x}u)^{2}\dot{u}. \label{Higer_estimates_2}%
\end{align}

\begin{remarka}
The last term of the above identity, $\frac{1}{2}\int_{\mathbb{R}}\partial
_{x}\rho^{\alpha}(\partial_{x}u)^{2}\dot{u}$ will be apear with sign minus in
the next identity
\end{remarka}

Let us observe that%
\begin{align}
& \int_{\mathbb{R}}(\partial_{x}\rho_{t}^{\gamma}+u\partial_{xx}^{2}%
\rho^{\gamma})\dot{u} =-\int_{\mathbb{R}}\rho_{t}^{\gamma}\partial_{x}\dot
{u}+\int_{\mathbb{R}}u\partial_{xx}^{2}\rho^{\gamma}\dot{u}\nonumber\\
&  =\int_{\mathbb{R}}u\partial_{x}\rho^{\gamma}\partial_{x}\dot{u}+\gamma
\int_{\mathbb{R}}\rho^{\gamma}\partial_{x}u\partial_{x}\dot{u}+\int
_{\mathbb{R}}u\partial_{xx}^{2}\rho^{\gamma}\dot{u}\nonumber\\
&  =-\int_{\mathbb{R}}\partial_{x}u\partial_{x}\rho^{\gamma}\dot{u}+\gamma
\int_{\mathbb{R}}\rho^{\gamma}\partial_{x}u\partial_{x}\dot{u}\nonumber\\
&  =\int_{\mathbb{R}}\partial_{x}u\rho\dot{u}^{2}-\int_{\mathbb{R}}%
\partial_{x}u\partial_{x}\left(  \rho^{\alpha}\partial_{x}u\right)  \dot
{u}+\gamma\int_{\mathbb{R}}\rho^{\gamma}\partial_{x}u\partial_{x}\dot
{u},\nonumber\\
&  =\int_{\mathbb{R}}\partial_{x}u\rho\dot{u}^{2}+\int_{\mathbb{R}}%
\rho^{\alpha}\partial_{x}u\partial_{x}\left(  \dot{u}\partial_{x}u\right)
+\gamma\int_{\mathbb{R}}\rho^{\gamma}\partial_{x}u\partial_{x}\dot
{u},\nonumber\\
&  =\int_{\mathbb{R}}\partial_{x}u\rho\dot{u}^{2}+\int_{\mathbb{R}}%
\rho^{\alpha}(\partial_{x}u)^{2}\partial_{x}\dot{u}+\int_{\mathbb{R}}\dot
{u}\rho^{\alpha}\partial_{x}u\partial_{xx}^{2}u+\gamma\int_{\mathbb{R}}%
\rho^{\gamma}\partial_{x}u\partial_{x}\dot{u},\nonumber\\
&  =\int_{\mathbb{R}}\partial_{x}u\rho\dot{u}^{2}+\int_{\mathbb{R}}%
\rho^{\alpha}(\partial_{x}u)^{2}\partial_{x}\dot{u}-\frac{1}{2}\int
_{\mathbb{R}}\partial_{x}(\dot{u}\rho^{\alpha})(\partial_{x}u)^{2}+\gamma
\int_{\mathbb{R}}\rho^{\gamma}\partial_{x}u\partial_{x}\dot{u},\nonumber\\
&  =\int_{\mathbb{R}}\partial_{x}u\rho\dot{u}^{2}+\frac{1}{2}\int_{\mathbb{R}%
}\rho^{\alpha}(\partial_{x}u)^{2}\partial_{x}\dot{u}-\frac{1}{2}%
\int_{\mathbb{R}}\dot{u}\partial_{x}\rho^{\alpha}(\partial_{x}u)^{2}%
+\gamma\int_{\mathbb{R}}\rho^{\gamma}\partial_{x}u\partial_{x}\dot{u},
\label{Higer_estimates_3}%
\end{align}
where we have used the equation of the velocity to replace
\[
-\partial_{x}\rho^{\gamma}=\rho\dot{u}-\partial_{x}(\rho^{\alpha}\partial
_{x}u).
\]
We sum up the relations $\left(  \text{\ref{Higer_estimates_1}}\right)  $,
$\left(  \text{\ref{Higer_estimates_2}}\right)  $ and $\left(
\text{\ref{Higer_estimates_3}}\right)  $ in order to obtain that:%
\[
\frac{1}{2}\frac{d}{dt}\int_{\mathbb{R}}\rho\dot{u}^{2}+\int_{\mathbb{R}}%
\rho^{\alpha}(\partial_{x}\dot{u})^{2}=\left(  \alpha+1\right)  \int
_{\mathbb{R}}\rho^{\alpha}(\partial_{x}u)^{2}\partial_{x}\dot{u}-\gamma
\int_{\mathbb{R}}\rho^{\gamma}\partial_{x}u\partial_{x}\dot{u}.
\]
Multiplying with $\sigma\left(  t\right)  $ and integrating in time on $[0,t]$
with $t>0$ leads to:%
\begin{align}
B\left(  \rho,u\right)  \left(  t\right)   &  =\frac{1}{2}\int_{\mathbb{R}%
}\sigma\left(  t\right)  \rho\dot{u}^{2}\left(  t\right)  +\int_{0}^{t}%
\int_{\mathbb{R}}\sigma\left(  t\right)  \rho^{\alpha}(\partial_{x}\dot
{u})^{2}\nonumber\\
&  =\int_{0}^{\min\left(  1,t\right)  }\int_{\mathbb{R}}\rho\dot{u}%
^{2}+\left(  \alpha+1\right)  \int_{0}^{t}\int_{\mathbb{R}}\sigma\rho^{\alpha
}(\partial_{x}u)^{2}\partial_{x}\dot{u}-\gamma\int_{0}^{t}\int_{\mathbb{R}%
}\sigma\rho^{\gamma}\partial_{x}u\partial_{x}\dot{u}. \label{A2}%
\end{align}
Obviously using $\left(  \text{\ref{estimate_for_A(rho,u)}}\right)  $ we have
that,
\begin{equation}
\int_{0}^{\min\left(  1,t\right)  }\int_{\mathbb{R}}\rho\dot{u}^{2}\leq
A\left(  \rho,u\right)  \left(  1\right)  \leq C. \label{estimate_A2_0}%
\end{equation}
for all $t>0$. Next, we infer using (\ref{Linf}) that:%
\begin{align}
\gamma\int_{0}^{t}\int_{\mathbb{R}}\sigma\rho^{\gamma}\partial_{x}%
u\partial_{x}\dot{u}  &  \leq\gamma\left\Vert \rho^{\gamma-\alpha}\right\Vert
_{L_{t}^{\infty}L^{\infty}}\left(  \int_{0}^{t}\int_{\mathbb{R}}\rho^{\alpha
}(\partial_{x}u)^{2}\right)  ^{\frac{1}{2}}\left(  \int_{0}^{t}\int
_{\mathbb{R}}\sigma^{2} \rho^{\alpha}(\partial_{x}\dot{u})^{2}\right)
^{\frac{1}{2}}\nonumber\\
&  \leq C\left(  t\right)  +\frac{1}{4}B\left(  \rho,u\right)  \left(
t\right)  , \label{estimate_A2_1}%
\end{align}
with $C$ a continuous increasing function. Finally, using again $\left(
\text{\ref{estimate_for_A(rho,u)}}\right)  $, (\ref{energy}) and
(\ref{loin_de_vide}), we get:%
\begin{align*}
\left(  \alpha+1\right)  \int_{0}^{t}\int_{\mathbb{R}}\sigma\rho^{\alpha
}(\partial_{x}u)^{2}\partial_{x}\dot{u}  &  \leq\frac{1}{4}\int_{0}^{t}%
\int_{\mathbb{R}}\sigma\rho^{\alpha}(\partial_{x}\dot{u})^{2}+\left(
\alpha+1\right)  ^{2}\int_{0}^{t}\int_{\mathbb{R}}\sigma\rho^{\alpha}%
(\partial_{x}u)^{4}\\
&  \leq\frac{1}{4}B\left(  \rho,u\right)  \left(  t\right)  +\left(
\alpha+1\right)  ^{2}\left\Vert \frac{1}{\rho}\right\Vert _{L^{\infty}%
_{t}(L^{\infty})}^{2\alpha} \int_{0}^{t}\int_{\mathbb{R}}\sigma\rho^{3\alpha
}(\partial_{x}u)^{4}\\
&  \leq\frac{1}{4}B\left(  \rho,u\right)  \left(  t\right)  +C\left(
t\right)  \int_{0}^{t}\sigma\left\Vert \rho^{\alpha}\partial_{x}u\right\Vert
_{L^{\infty}}^{2}\int_{\mathbb{R}}\rho^{\alpha}(\partial_{x}u)^{2}\\
&  \leq\frac{1}{4}B\left(  \rho,u\right)  \left(  t\right)  +C\left(
t\right)  \sup_{\tau\in\left[  0,t\right]  }\sigma\left(  \tau\right)
\left\Vert (\rho^{\alpha}\partial_{x}u)\left(  \tau\right)  \right\Vert
_{L^{\infty}}^{2}.
\end{align*}
Let us observe that for all $t>0$ we have using Gagliardo-Nirenberg
inequality, (\ref{Linf}) and (\ref{loin_de_vide}):%
\begin{align}
\sigma\left(  t\right)  \left\Vert \rho^{\alpha}\partial_{x}u\left(  t\right)
\right\Vert _{L^{\infty}}^{2}  &  \leq2\sigma\left\Vert \left(  \rho^{\alpha
}\partial_{x}u-(\rho^{\gamma}-1)\right)  \left(  t\right)  \right\Vert
_{L^{\infty}}^{2}+2\left\Vert \left(  \rho^{\gamma}-1\right)  \left(
t\right)  \right\Vert _{L^{\infty}}^{2}\nonumber\\
&  \leq2\sigma\left\Vert \left(  \rho^{\alpha}\partial_{x}u-(\rho^{\gamma
}-1)\right)  \left(  t\right)  \right\Vert _{L^{2}}\left\Vert \partial
_{x}\left(  \rho^{\alpha}\partial_{x}u-(\rho^{\gamma}-1)\right)  \left(
t\right)  \right\Vert _{L^{2}}+C\left(  t\right) \nonumber\\
&  \leq2\sigma\left(  \left\Vert \rho^{\alpha}\partial_{x}u\right\Vert
_{L^{2}}+C\left(  t\right)  \left\Vert \rho-1\right\Vert _{L^{2}}\right)
\left\Vert \rho\dot{u}\right\Vert _{L^{2}}+C\left(  t\right) \nonumber\\
&  \leq C\left(  t\right)  \left(  \sigma^{\frac{1}{2}}\left\Vert \rho
^{\alpha}\partial_{x}u\right\Vert _{L^{2}}+C\left(  t\right)  \right)
\sigma^{\frac{1}{2}}\left\Vert \rho^{\frac{1}{2}}\dot{u}\right\Vert _{L^{2}%
}+C\left(  t\right) \nonumber\\
&  \leq C\left(  t\right)  \left(  A^{\frac{1}{2}}\left(  \rho,u\right)
\left(  t\right)  +C\left(  t\right)  \right)  B^{\frac{1}{2}}\left(
\rho,u\right)  \left(  t\right)  +C\left(  t\right) .
\label{flux_eff_LinfLinf}%
\end{align}
Thus, we get from (\ref{estimate_for_A(rho,u)}) and Young inequality:
\begin{gather}
\left(  \alpha+1\right)  \int_{0}^{t}\int_{\mathbb{R}}\sigma\rho^{\alpha
}(\partial_{x}u)^{2}\partial_{x}\dot{u}\leq\frac{1}{4}B\left(  \rho,u\right)
\left(  t\right)  +C\left(  t\right)  \left(  A^{\frac{1}{2}}\left(
\rho,u\right)  \left(  t\right)  +C\left(  t\right)  \right)  B^{\frac{1}{2}%
}\left(  \rho,u\right)  \left(  t\right)  +C\left(  t\right) \nonumber\\
\leq C\left(  t\right)  +\frac{1}{2}B\left(  \rho,u\right)  \left(  t\right) .
\label{estimate_A2_2}%
\end{gather}
Gathering $\left(  \text{\ref{estimate_A2_0}}\right)  $, $\left(
\text{\ref{estimate_A2_1}}\right)  $ and $\left(  \text{\ref{estimate_A2_2}%
}\right)  $ yields the fact that $B$ is also bounded:%
\begin{equation}
B\left(  \rho,u\right)  \left(  t\right)  \leq C\left(  t\right) ,
\label{estimate_for_A2}%
\end{equation}
with $C$ a continuous increasing function. The control over $\left\Vert
\frac{1}{\rho}\right\Vert _{L^{\infty}}$, $A\left(  \rho,u\right)  $ and
$B\left(  \rho,u\right)  $ gives us, via the estimate $\left(
\text{\ref{flux_eff_LinfLinf}}\right)  $ the following
\begin{equation}
\sigma\left(  t\right)  ^{\frac{1}{2}}\left\Vert \partial_{x}u(t)\right\Vert
_{L^{\infty}}\leq C\left(  t\right)  , \label{u_x_in_linfty}%
\end{equation}
for any $t\geq0$.

\subsection{Uniform BV-estimates for the effective velocities $v_{n}$}

Owing to the estimate $\left(  \text{\ref{loin_de_vide}}\right)  $ and
(\ref{energy}) we recover the following estimates:%
\[
\left\Vert \partial_{x}u_{n}\right\Vert _{L_{t}^{2}(L^{2})}\leq\left\Vert
\frac{1}{\rho_{n}}\right\Vert _{L_{t}^{\infty}(L^{\infty})}^{\frac{\alpha}{2}%
}\leq C(t),\;\;\left\Vert \sqrt{\rho_{n}}u_{n}\right\Vert _{L_{t}^{2}(L^{2}%
)}\leq C_{1}\left(  t\right)  ,
\]
where $C$, $C_{1}$ are increasing continuous functions. From Sobolev
embedding, we get that for any $t>0$, there exists $C\left(  t\right)  $ such
that
\begin{equation}
\left\Vert u_{n}\right\Vert _{L_{t}^{1}(L^{\infty})}\leq C\left(  t\right)
.\label{u_borne}%
\end{equation}
Let us introduce the flow of $u_{n}$ i.e.%
\begin{equation}
X_{n}\left(  t,x\right)  =x+\int_{0}^{t}u_{n}\left(  \tau,X_{n}\left(
\tau,x\right)  \right)  d\tau.\label{flow}%
\end{equation}
We immediately get that:
\[
-\left\vert x\right\vert -C\left(  t\right)  \leq\left\vert X_{n}^{\pm
1}\left(  t,x\right)  \right\vert \leq\left\vert x\right\vert +C\left(
t\right)  ,
\]
which implies that for any $L>0$ the segment
\[
\left[  X_{n}^{-1}\left(  t,-L\right)  ,X_{n}^{-1}\left(  t,L\right)  \right]
\leq\left[  -L-C\left(  t\right)  ,L+C\left(  t\right)  \right]  .
\]
This information is usefull in order to show that we can propagate the
$L_{loc}^{\infty}$ norm of $v_{n}$. Indeed, let us recall that:
\[
\partial_{t}v_{n}+u_{n}\partial_{x}v_{n}+\frac{P^{\prime}\left(  \rho
_{n}\right)  \rho_{n}^{2}}{\mu\left(  \rho_{n}\right)  }\frac{\mu\left(
\rho_{n}\right)  }{\rho_{n}^{2}}\partial_{x}\rho_{n}=0,
\]
rewrites as%
\begin{equation}
\partial_{t}v_{n}+u_{n}\partial_{x}v_{n}+\frac{P^{\prime}\left(  \rho
_{n}\right)  \rho_{n}^{2}}{\mu\left(  \rho_{n}\right)  }\left(  v_{n}%
-u_{n}\right)  =0.\label{v_eq}%
\end{equation}
Passing into Lagrangian coordinates (see Appendix) i.e.
\[
\left(  \widetilde{v}_{n},\widetilde{u}_{n},\widetilde{\rho}_{n}\right)
\left(  t,x\right)  =\left(  v_{n},u_{n},\rho_{n}\right)  \left(
t,X_{n}\left(  t,x\right)  \right)  ,
\]
we see that $\left(  \text{\ref{v_eq}}\right)  $ rewrites as:%
\begin{equation}
\partial_{t}\widetilde{v}_{n}+\frac{P^{\prime}\left(  \widetilde{\rho}%
_{n}\right)  \widetilde{\rho_{n}}^{2}}{\mu\left(  \widetilde{\rho}_{n}\right)
}\widetilde{v}_{n}=\frac{P^{\prime}\left(  \widetilde{\rho}_{n}\right)
\widetilde{\rho_{n}}^{2}}{\mu\left(  \widetilde{\rho}_{n}\right)  }%
\widetilde{u}_{n}.\label{v_in_lagrange}%
\end{equation}
The last relation implies using (\ref{Linf}) and (\ref{loin_de_vide}):%
\begin{align*}
|\widetilde{v}_{n}\left(  t,x\right)  | &  \leq\left\vert v_{0n}\left(
x\right)  \exp\left(  -\int_{0}^{t}\frac{P^{\prime}\left(  \widetilde{\rho
}_{n}\left(  \tau,x\right)  \right)  \widetilde{\rho}_{n}^{2}(\tau,x)}%
{\mu\left(  \widetilde{\rho}_{n}\left(  \tau,x\right)  \right)  }d\tau\right)
\right\vert \\
&  +\left\vert \int_{0}^{t}\exp\left(  -\int_{s}^{t}\frac{P^{\prime}\left(
\widetilde{\rho}_{n}\left(  \tau,x\right)  \widetilde{\rho}_{n}^{2}%
(\tau,x)\right)  }{\mu\left(  \widetilde{\rho}_{n}\left(  \tau,x\right)
\right)  }d\tau\right)  \frac{P^{\prime}\left(  \widetilde{\rho}_{n}\left(
s,x\right)  \right)  \widetilde{\rho}_{n}^{2}(s,x)}{\mu\left(  \widetilde
{\rho}_{n}\left(  s,x\right)  \right)  }\widetilde{u}_{n}\left(  s,x\right)
ds\right\vert \\
&  \leq C\left(  t\right)  \left(  |v_{0n}\left(  x\right)  |+\int_{0}%
^{t}\left\Vert u_{n}\left(  s\right)  \right\Vert _{L^{\infty}}ds\right)  .\\
&  \leq C\left(  t\right)  \left(  1+|v_{0n}\left(  x\right)  |\right)  ,
\end{align*}
and consequently for any $t>0$, $x\in\mathbb{R}$:
\[
|v_{n}\left(  t,x\right)  |\leq C\left(  t\right)  \left(  1+|v_{0n}\left(
X_{n}^{-1}\left(  t,x\right)  \right)  |\right)  .
\]
Thus, we see that:%
\begin{equation}
\left\Vert v_{n}\left(  t\right)  \right\Vert _{L^{\infty}\left(  \left[
-L,L\right]  \right)  }\leq C\left(  t\right)  \left(  1+\left\Vert
v_{0n}\right\Vert _{L^{\infty}\left(  \left[  -L-C\left(  t\right)
,L+C\left(  t\right)  \right]  \right)  }\right)  .\label{v_local_borne}%
\end{equation}
In addition $(v_{0}^{n})_{n\in\mathbb{N}}$ is uniformly bounded in
$L_{loc}^{\infty}(\mathbb{R})$. Indeed since $v_{0}$ is in $L_{loc}^{\infty
}(\mathbb{R})$ (see the Remark \ref{remimpo}), we have for any $x\in\lbrack-L,L]$ and any $n\in\mathbb{N}$:
\begin{equation}
|v_{0}^{n}(x)|\leq\int_{-1}^{1}j(y)v_{0}(x-\frac{y}{n})dy\leq\Vert v_{0}%
\Vert_{L^{\infty}([-L-1,L+1])}.\label{3.60}%
\end{equation}
This piece of information along with the estimate:%
\[
\partial_{x}v_{n}\left(  t,x\right)  \leq C\left(  t\right)
\]
ensures that $v^{n}$ is uniformly bounded in $L^{\infty}([0,T];BV_{loc}%
(\mathbb{R}))$. Ideed, the function
\[
w_{n}\left(  t,x\right)  =v_{n}\left(  t,x\right)  -C\left(  t\right)  x
\]
being nonincreasing, it holds using (\ref{v_local_borne}) that:%
\[
TV_{\left[  -L,L\right]  }w_{n}(t,\cdot)=v_{n}\left(  t,-L\right)
-v_{n}\left(  t,L\right)  +2C\left(  t\right)  L\leq C\left(  t\right)
\left(  L+\left\Vert v_{0n}\right\Vert _{L^{\infty}\left(  \left[  -L-C\left(
t\right)  ,L+C\left(  t\right)  \right]  \right)  }\right)  .
\]
Owing to the fact that%
\[
v_{n}=w_{n}\left(  t,x\right)  +C\left(  t\right)  x
\]
we get that:%
\begin{equation}
TV_{\left[  -L,L\right]  }v_{n}(t,\cdot)\leq C\left(  t\right)  \left(
L+\left\Vert v_{0n}\right\Vert _{L^{\infty}\left(  \left[  -L-C\left(
t\right)  ,L+C\left(  t\right)  \right]  \right)  }\right)  .\label{BV_en_x}%
\end{equation}
From $\left(  \text{\ref{v_local_borne}}\right)  $, (\ref{3.60}),
(\ref{BV_en_x}), we get :
\begin{equation}
\left\Vert v_{n}\left(  t\right)  \right\Vert _{BV\left(  \left[  -L,L\right]
\right)  }\leq C\left(  T,L\right)  .\label{vn_bv(t)}%
\end{equation}
Owing to (\ref{Linf}), (\ref{loin_de_vide}), $\left(  \text{\ref{u_borne}%
}\right)  $, $\left(  \text{\ref{v_local_borne}}\right)  $ and $\left(
\text{\ref{v_in_lagrange}}\right)  $ we get that:%
\begin{equation}
\left\Vert \partial_{t}\widetilde{v}_{n}\right\Vert _{L^{1}\left(
[0,T]\times\left[  -L,L\right]  \right)  }\leq C\left(  T,L\right)
.\label{vn_lag_deriv_bornee}%
\end{equation}
Next, fix $\phi\in C_{b}\left(  \left[  0,T\right]  \times\left[  -L,L\right]
\right)  $ with:
\[
\left\Vert \phi\right\Vert _{L^{\infty}\left(  \left[  0,T\right]
\times\left[  -L,L\right]  \right)  }\leq1
\]
and let us write that:%
\begin{align*}
&  \int_{0}^{T}\int_{-L}^{L}\phi\left(  t,x\right)  \partial_{t}v_{n}\left(
t,x\right)  dxdt\\
&  =\int_{0}^{T}\left(  \int_{X_{n}(t,-L)}^{X_{n}(t,L)}\widetilde{\phi}\left(
t,x\right)  \widetilde{\partial_{t}v}_{n}\left(  t,x\right)  \frac
{\widetilde{\rho}_{n}\left(  t,x\right)  }{\rho_{0n\left(  x\right)  }%
}dx\right)  dt\\
&  =\int_{0}^{T}\left(  \int_{X_{n}(t,-L)}^{X_{n}(t,L)}\widetilde{\phi}\left(
t,x\right)  \partial_{t}\widetilde{v}_{n}\left(  t,x\right)  \frac
{\widetilde{\rho}_{n}\left(  t,x\right)  }{\rho_{0n\left(  x\right)  }%
}dx\right)  dt\\
&  -\int_{0}^{T}\left(  \int_{X_{n}(t,-L)}^{X_{n}(t,L)}\widetilde{\phi}\left(
t,x\right)  \widetilde{u}_{n}\left(  t,x\right)  \widetilde{\partial_{x}v_{n}%
}\left(  \left(  t,x\right)  \right)  \frac{\widetilde{\rho}_{n}\left(
t,x\right)  }{\rho_{0n\left(  x\right)  }}dx\right)  dt\\
&  =\int_{0}^{T}\left(  \int_{X_{n}(t,-L)}^{X_{n}(t,L)}\widetilde{\phi}\left(
t,x\right)  \partial_{t}\widetilde{v}_{n}\left(  t,x\right)  \frac
{\widetilde{\rho}_{n}\left(  t,x\right)  }{\rho_{0n\left(  x\right)  }%
}dx\right)  dt-\int_{0}^{T}\int_{-L}^{L}\phi\left(  t,x\right)  u_{n}\left(
t,x\right)  \partial_{x}v_{n}\left(  \left(  t,x\right)  \right)  dxdt.
\end{align*}
Owing to (\ref{Linf}), (\ref{loin_de_vide}), $\left(  \text{\ref{u_borne}%
}\right)  $, $\left(  \text{\ref{vn_bv(t)}}\right)  $ and $\left(
\text{\ref{vn_lag_deriv_bornee}}\right)  $ and using the fact that $\phi
u_{n}$ belongs to $L^{1}([0,T],C^{0}(\mathbb{R}))$, we conclude that
\begin{equation}
\left\vert \int_{0}^{T}\int_{-L}^{L}\phi\left(  t,x\right)  \partial_{t}%
v_{n}\left(  t,x\right)  dxdt\right\vert \leq C\left(  T,L\right)
.\label{r34}%
\end{equation}
Combining (\ref{r34}) and $\left(  \text{\ref{BV_en_x}}\right)  $ gives us for
any $T>0$, $L>0$:%
\begin{equation}
\left\Vert v_{n}\right\Vert _{BV\left(  [0,T]\times\left[  -L,L\right]
\right)  }\leq C\left(  T,L\right)  .\label{BV1d}%
\end{equation}

\subsection{Compactness}

We recall the previous estimates that we have obtained, for every $T>0$ we
have for $C$ a continuous increasing function independent on $n$ and any
$n\in\mathbb{N}$:
\begin{equation}
C(T)^{-1}\leq\rho_{n}(T,\cdot)\leq C(T),\label{n1.17}%
\end{equation}%
\begin{equation}
\begin{aligned} &\sup_{0<t\leq T}\big( \|\rho_n(t,\cdot)-1\|_{L^2}+\|u_n(t,\cdot)\|_{L^2}+\|\partial_x\rho_n(t,\cdot)\|_{L^2}+\sigma(t)^{\frac{1}{2}}\|\partial_x u_n(t,\cdot)\|_{L^2}\\ &+\sigma(t)^{\frac{1}{2}}(\|\dot{u_n}(t,\cdot)\|_{L^2}+\|\partial_x(\rho_n^\alpha\partial_x u_n(t,\cdot)-P(\rho_n)+P(1))\|_{L^2}\big)\leq C(T), \end{aligned}\label{n1.18}%
\end{equation}%
\begin{equation}
\begin{aligned} &\int^T_0 [\|\partial_x u_n(t,\cdot)\|_{L^2}^2+\|\partial_x\rho_n(t,\cdot)\|_{L^2}^2+\sigma(t)\|\dot{u_n}(t,\cdot)\|_{L^2}^2+\sigma(t)\|\partial_x \dot{u_n}(t,\cdot)\|_{L^2}^2]dt \leq C(T), \end{aligned}\label{1.19n}%
\end{equation}%
\begin{equation}
\int_{0}^{T}\sigma^{\frac{1}{2}}\left(  \tau\right)  \left\Vert \partial
_{x}u_{n}\left(  \tau\right)  \right\Vert _{L^{\infty}}^{2}d\tau\leq C\left(
T\right)  .\label{nLipo}%
\end{equation}%
\begin{equation}
\sup_{0<t\leq T}\sigma(t)^{\frac{1}{2}}\Vert\partial_{x}u_{n}(t,\cdot
)\Vert_{L^{\infty}}\leq C\left(  T\right)  .\label{nLipo1}%
\end{equation}
Using classical arguments (see \cite{Jiu,MV}), we prove that up to a
subsequence, $(\rho_{n},u_{n})_{n\in\mathbb{N}}$ converges in the sense of
distributions to $(\rho,u)$, a global weak solution of (\ref{0.1}).
Furthermore the limit functions $\rho$, $u$ inherit all the bounds
(\ref{n1.17}), (\ref{n1.18}), (\ref{1.19n}), (\ref{nLipo}), (\ref{BV2}) and
(\ref{nLipo1}) via Fatou type-lemmas for the weak topology.
\\
We wish now to prove (\ref{uBV}), to do this we are going to prove that up to a subsequence $(v_n)_{n\in\mathbb{N}}$ converges almost everywhere to $v$ on $\R^+\times\R$. This is a direct consequence of the estimate (\ref{BV1d}), indeed since $(v_n)_{n\in\mathbb{N}}$ is uniformly bounded in $BV_{loc}((0,T)\times\R)$ for any $T>0$, we deduce that up to a subsequence $(v_n)_{n\in\mathbb{N}}$ converges to $v$ in $L^1_{loc}((0,T)\times\R)$. In particular up to a subsequence $(v_n)_{n\in\mathbb{N}}$ converges almost everywhere to $v$ in $(0,+\infty)\times\R.$
Using now (\ref{pointfin1}) and the fact that $v_n$  converges almost everywhere to $v$ on $\R^+\times\R$ implies (\ref{uBV}) since for all $x>y$ and $t>0$ we have:
$$
\frac{v_n(t,x)-v_n(t,y)}{x-y}=\frac{1}{x-y}\int^x_y\p_z v_n(t,z) dz\leq C(t),
$$
with $C$ a continuous function on $\R^+$.  It concludes the proof of (\ref{uBV}).
\subsection{Uniqueness}

Consider two solutions $\left(  \rho_{i},u_{i}\right)  $ ( $i\in\overline
{1,2}$) verifying the estimates (\ref{1.17})-(\ref{Lipo1}) and generated by
the same initial data:
\begin{equation}
\left\{
\begin{array}
[c]{l}%
\partial_{t}\rho_{i}+\partial_{x}\left(  \rho_{i}u_{i}\right)  =0,\\
\partial_{t}\left(  \rho_{i}u_{i}\right)  +\partial_{x}\left(  \rho_{i}%
u_{i}^{2}\right)  -\partial_{x}(\mu(\rho_{i})\partial_{x}u_{i})+\partial
_{x}p_{i}=0,\\
\left(  \rho_{i|t=0},u_{i|t=0}\right)  =\left(  \rho_{0},u_{0}\right)  .
\end{array}
\right.  \label{NSCM_i}%
\end{equation}
We define now the flows generated by $u_{i}$
\[
X_{i}(t,x)=x+\int_{0}^{t}u_{i}\left(  \tau,X\left(  \tau,x\right)  \right)
d\tau
\]
and denoting with tildes the functions%
\[
\widetilde{v}_{i}\left(  t,x\right)  =v_{i}\left(  t,X_{i}\left(  t,x\right)
\right)
\]
for $v\in\left\{  \rho,u\right\}  $. We get that (according to the results
from the Appendix):%
\begin{equation}
\left\{
\begin{array}
[c]{l}%
\partial_{t}\left(  \dfrac{\partial X^{i}}{\partial x}\widetilde{\rho}%
_{i}\right)  =0,\\
\rho_{0}\partial_{t}\widetilde{u}_{i}-\partial_{x}\left(  \dfrac
{\widetilde{\rho}_{i}\mu(\widetilde{\rho}_{i})}{\rho_{0}}\,\partial
_{x}\widetilde{u}_{i}\right)  +\partial_{x}P(\widetilde{\rho}_{i})=0,\\
X_{i}(t,x)=x+%
{\displaystyle\int_{0}^{t}}
\widetilde{u}_{i}\left(  \tau,x\right)  d\tau.
\end{array}
\right.  \label{NSCM_lag_i}%
\end{equation}
for $i=1,2$. Setting $\delta\widetilde{u}=\widetilde{u}_{1}-\widetilde{u}_{2}%
$, by difference we have that:%
\begin{equation}
\rho_{0}\partial_{t}\delta\widetilde{u}-\partial\left(  \dfrac{\widetilde
{\rho}_{1}\mu(\widetilde{\rho}_{1})}{\rho_{0}}\partial_{x}\delta\widetilde
{u}\right)  =\partial_{x}G_{1}+\partial_{x}G_{2},\label{difference}%
\end{equation}
where%
\[
\left\{
\begin{array}
[c]{l}%
G_{1}=P\left(  \dfrac{\rho_{0}}{1+\int_{0}^{t}\partial_{x}[\widetilde{u}_{2}%
]}\right)  -P\left(  \dfrac{\rho_{0}}{1+\int_{0}^{t}\partial_{x}[\widetilde
{u}_{1}]}\right)  ,\\
G_{2}=\left(  \dfrac{\widetilde{\rho}_{1}\mu(\widetilde{\rho}_{1})}{\rho_{0}%
}-\dfrac{\widetilde{\rho}_{2}\mu(\widetilde{\rho}_{2})}{\rho_{0}}\right)
\partial_{x}\widetilde{u}_{2}.
\end{array}
\right.
\]
We multiply $\left(  \text{\ref{difference}}\right)  $ by $\delta\widetilde
{u}$, integrate it over $\mathbb{R}$ and by obvious manipulation we get for
$t>0$:
\begin{equation}
\begin{aligned} &\frac{1}{2}\int_{\mathbb{R}}\rho_{0}(x)(\delta\widetilde{u})^{2}(t,x) dx+\frac{1}{2}\inf_{s\in]0,t],x}\dfrac{\tilde{\rho}_{1}\left( s,x\right) \mu(\tilde{\rho}_{1}\left( s,x\right) )}{\rho_{0}\left( x\right) }\int_{0}^{t}\int_{\mathbb{R}}\left( \partial_{x}(\delta\widetilde{u}(s,x))\right) ^{2}ds\, dx\\ &\hspace{6cm}\leq C(t)[\int_{0}^{t}\int_{\mathbb{R}}(G_{1})^{2}+\int_{0}^{t}\int_{\mathbb{R}}(G_{2})^{2}], \end{aligned}\label{uni1}%
\end{equation}
with $C$ a continuous increasing function. In the following we will estimate
$G_{1}$ and $G_{2}$. First, we we get using (\ref{flow_ineg2}) for $t>0$ and
$x\in\mathbb{R}$:%
\begin{equation}
\begin{aligned} \delta\widetilde{\rho}\left( t,x\right) &=\delta{\rho}_1-\delta{\rho}_2=\dfrac{\rho_{0}\left( x\right) }{1+\int_{0}^{t}\partial_{x}\tilde{u}_{1}}-\dfrac{\rho_{0}\left( x\right) }{1+\int_{0}^{t}\partial_{x}\tilde{u}_{2}}=\frac{-\rho_{0}\left( x\right) \int_{0}^{t}\partial_{x}\delta\tilde{u}\left( \tau,x\right) d\tau}{\p_x X^{1}(t,x)\p_x X^{2}(t,x)}\\ |\delta\widetilde{\rho}\left( t,x\right) |&\leq\sqrt{t}C\left( t\right) \left( \int_{0}^{t}\left\vert \partial_{x}\delta\widetilde{u}\left( \tau,x\right) \right\vert ^{2}d\tau\right) ^{\frac{1}{2}}\end{aligned}\label{3.54}%
\end{equation}
and consequently using (\ref{dens}) we get:%
\[
G_{1}\left(  t,x\right)  \leq\sup_{s\in\left[  1/C\left(  T\right)  ,C\left(
T\right)  \right]  }P^{\prime}\left(  s\right)  \frac{\rho_{0}\left(
x\right)  \left\vert \int_{0}^{t}\partial_{x}\delta\widetilde{u}\left(
\tau,x\right)  d\tau\right\vert }{|\partial_{x}X^{1}(t,x)\partial_{x}%
X^{2}(t,x)|}\leq\sqrt{t}C\left(  t\right)  \left(  \int_{0}^{t}\left\vert
\partial_{x}\delta\widetilde{u}\left(  \tau,x\right)  \right\vert ^{2}%
d\tau\right)  ^{\frac{1}{2}},
\]
with $C$ a continuous increasing function. It implies that
\begin{equation}
\int_{0}^{t}\int_{\mathbb{R}}(G_{1})^{2}\left(  s,x\right)  ds\,dx\leq
t^{\frac{3}{2}}C\left(  t\right)  \int_{0}^{t}\int_{\mathbb{R}}(\partial
_{x}\delta\widetilde{u})^{2}(s,x)dsdx.\label{ineg_G1}%
\end{equation}
Let us turn our attention towards $G_{2}$. We first write that for any
$(t,x)\in\mathbb{R}^{+}\times\mathbb{R}$, we have:%
\[
\begin{aligned}
G_{2}(t,x)&=\left(  \dfrac{\tilde{\rho}_{1}\mu(\tilde{\rho}_{1})}{\rho_{0}}%
-\dfrac{\tilde{\rho}_{2}\mu(\tilde{\rho}_{2})}{\rho_{0}}\right)  \partial
_{x}\tilde{u}_{2}(t,x)\\
&=\frac{1}{\rho_{0}(x)}\left(\mu\right)  ^{\prime}\left(
\theta_{t,x}\tilde{\rho}_{1}(t,x)+\left(  1-\theta_{t,x}\right)  \tilde{\rho}_{2}(t,x)\right)
\delta\tilde{\rho}(t,x)\partial_{x}\tilde{u}_{2}(t,x).%
\end{aligned}
\]
Thus, we get using (\ref{1.17}), (\ref{3.54}) that for $t>0$:%
\begin{align*}
& \left\vert G_{2}\left(  t,x\right)  \right\vert ^{2}\leq C\left(  t\right)
\left(  t\left\vert \partial_{x}\widetilde{u}_{2}\left(  t,x\right)
\right\vert ^{2}\right)  \left(  \int_{0}^{t}\left\vert \partial_{x}%
\delta\widetilde{u}\left(  \tau,x\right)  \right\vert ^{2}d\tau\right)  \\
& \leq C\left(  t\right)  \left(  \left(  \sigma^{\frac{1}{2}}\left(
t\right)  \mathbf{1}_{\left[  0,1\right]  }\left(  t\right)  +t\mathbf{1}%
_{[1,\infty)}\left(  t\right)  \right)  \left\Vert \partial_{x}\widetilde
{u}_{2}(t)\right\Vert _{L^{\infty}}\right)  ^{2}\left(  \int_{0}^{t}\left\vert
\partial_{x}\delta\widetilde{u}\left(  \tau,x\right)  \right\vert ^{2}%
d\tau\right)
\end{align*}
such that by integration and using (\ref{Lipo1}), (\ref{flow_ineg1}) and
(\ref{1.17}) we have:
\begin{equation}
\int_{\mathbb{R}}\left\vert G_{2}\left(  t,x\right)  dx\right\vert ^{2}\leq
C\left(  t\right)  \int_{0}^{t}\int_{\mathbb{R}}[\partial_{x}(\delta
\widetilde{u}(s,x))|^{2}dx\,ds.\label{uni2}%
\end{equation}
Putting together the inequalities (\ref{uni1}), $\left(  \text{\ref{ineg_G1}%
}\right)  $, (\ref{uni2}) and integrating in time, we get that for $t>0$:
\begin{equation}
\begin{aligned} &\frac{1}{2}\int_{\mathbb{R}}\rho_{0}(x)(\delta\widetilde{u})^{2}(t,x) dx+\frac{1}{2}\inf_{s\in]0,t],x}\dfrac{\tilde{\rho}_{1}\left( s,x\right) \mu(\tilde{\rho}_{1}\left( s,x\right) )}{\rho_{0}\left( x\right) }\int_{0}^{t}\int_{\mathbb{R}}\left( \partial_{x}(\delta\widetilde{u}(s,x))\right) ^{2}ds\, dx\\ &\hspace{6cm}\leq tC_1(t) \int_{0}^{t}\int_{\mathbb{R}}[\partial_x(\delta\widetilde{u}(s,x))|^{2}dx\,ds, \end{aligned}\label{uni3}%
\end{equation}
with $C_{1}$ a continuous increasing function. Taking $T_{0}>0$ small enough,
we have using a bootstrap argument for any $t\in\lbrack0,T_{0}]$:%
\[
\frac{1}{2}\int_{0}^{1}\rho_{0}(\delta\widetilde{u})^{2}+\frac{1}{4}\inf
_{t,x}\dfrac{\widetilde{\rho}_{1}\left(  t,x\right)  \mu(\widetilde{\rho}%
_{1}\left(  t,x\right)  )}{\rho_{0}\left(  x\right)  }\int_{0}^{t}\int_{0}%
^{1}\left(  \partial_{x}(\delta\widetilde{u})\right)  ^{2}\leq0\text{ }\forall
t\in\left[  0,T_{0}\right]  .
\]
Thus, we get a local uniqueness property. Reiterating this process gives us
the uniqueness of the two solutions on their whole domain of definition.

\section{Appendix}

In this appendix, we gather a few useful facts regarding the $1D$
Navier-Stokes equations in Lagrangian coordinates. The results belong to the
mathematical folklore and can be found in, by now classical papers devoted to
the $1D$ Navier-Stokes system, see \cite{KS77}, \cite{Ser86a}, \cite{Ser86b}.
The Lagrangian framework offers an elegant method of obtaining apriori
estimates (for example on the $L^{\infty}$ norm of $\rho$) either uniqueness
of solutions (see the relatively recent paper  \cite{Fourier}). \newline\ Let
us first derive the Lagrangian formulation of the Navier-Stokes system. We
will supose that we are give $\left(  \rho,u\right)  \in L^{\infty}\left(
[0,\infty)\times\mathbb{R}\right)  \times L^{\infty}\left(  L^{2}%
(\mathbb{R})\right)  \cap L^{2}\left(  \dot{H}\left(  \mathbb{R}\right)
\right)  $ a solution of the Navier-Stokes system
\begin{equation}
\left\{
\begin{array}
[c]{l}%
\rho_{t}+\partial_{x}\left(  \rho u\right)  =0,\\
\partial_{t}\left(  \rho u\right)  +\partial_{x}\left(  \rho u^{2}\right)
-\partial_{x}\left(  \mu\left(  \rho\right)  \partial_{x}u\right)
+\partial_{x}P\left(  \rho\right)  =0,\\
\left(  \rho_{|t=0},u_{|t=0}\right)  =\left(  \rho_{0},u_{0}\right)  .
\end{array}
\right.  \label{Sistem_appendix}%
\end{equation}
First, we recall the definition of the flow of $u$. 

\begin{proposition}
Consider $T>0$ and $u\in L^{2}\left(  (0,T);L^{\infty}\left(  \mathbb{R}%
\right)  \right)  $ with $\partial_{x}u\in L^{1}\left(  (0,T);L^{\infty
}\left(  \mathbb{R}\right)  \right)  $. Then, for any $x\in\mathbb{R}$ there
exists a unique solution $X\left(  \cdot,x\right)  :[0,\infty)\rightarrow
\mathbb{R}$ of%
\begin{equation}
\left\{
\begin{array}
[c]{l}%
X\left(  t,x\right)  =x+%
{\displaystyle\int_{0}^{t}}
u\left(  t,X\left(  t,x\right)  \right)  ,\\
X\left(  0,x\right)  =x.
\end{array}
\right.  \label{FLOW}%
\end{equation}
Moreover $X\left(  t,x\right)  $ verifies the following properties:

\begin{itemize}
\item $X\in BV_{loc}\left(  \left[  0,T\right]  \times\mathbb{R}\right)  $ for
any $T>0$. In addition, for all $t\geq0$ and for almost all $x\in\mathbb{R}$
\[
\partial_{x}X\left(  t,x\right)  =\exp\left(  \int_{0}^{t}\partial_{x}u\left(
\tau,X\left(  \tau,x\right)  \right)  d\tau\right)
\]

\item For each $t>0$, $X\left(  t,\cdot\right)  $ is a homeorphism from
$\mathbb{R}$ to $\mathbb{R}$.

\item We have that $\partial_{t}X,\partial_{t}X\in L_{t}^{2}(L_{x}^{\infty})$
and $\partial_{x}X,\partial_{x}X^{-1}\in L_{t}^{\infty}(L_{x}^{\infty})$
\end{itemize}
\end{proposition}

\begin{notation}
For any function $v:\left[  0,\infty\right)  \times\mathbb{R\rightarrow R}$,
we denote by $\widetilde{v}$ the function defined as:%
\[
\widetilde{v}\left(  t,x\right)  :=v\left(  t,X\left(  t,x\right)  \right)
\]

\end{notation}

We note that
\[
X\left(  t,x\right)  =x+\int_{0}^{t}u\left(  \tau,X(\tau,x)\right)
=x+\int_{0}^{t}\widetilde{u}\left(  \tau,x\right)  d\tau.
\]
and thus
\[
\frac{\partial X}{\partial x}\left(  t,x\right)  =1+\int_{0}^{t}\partial
_{x}\widetilde{u}\left(  \tau,x\right)  d\tau.
\]
Owing to%
\begin{equation}
\partial_{x}\widetilde{v}\left(  t,x\right)  =\widetilde{\partial_{x}v}\left(
t,x\right)  \frac{\partial X}{\partial x}\left(  t,x\right)
\label{derivative_lag_1}%
\end{equation}
we obtain that%
\begin{equation}
\widetilde{\partial_{x}v}\left(  t,x\right)  =\frac{\partial X}{\partial
x}\left(  t,x\right)  ^{-1}\partial_{x}\widetilde{v}\left(  t,x\right)
=\frac{1}{1+\int_{0}^{t}\partial_{x}\widetilde{u}\left(  \tau,x\right)  d\tau
}\partial_{x}\widetilde{v}\left(  t,x\right)  .\label{derivative_lag_2}%
\end{equation}

Let us investigate the first equation of $\left(  \text{\ref{Sistem_appendix}%
}\right)  $. For any $\psi\in\mathcal{D}\left(  (0,T)\times\mathbb{R}\right)
$ we have that :
\[
\int_{0}^{T}\int_{\mathbb{R}}\rho\psi_{t}+\rho u\partial_{x}\psi=0.
\]
Owing to the fact that $\rho$, $\rho u\in L_{T}^{2}\left(  L_{loc}^{2}\right)
$ the set of test functions can be enlarged to $\psi\in C^{0}\left(  \left(
0,T\right)  \times\mathbb{R}\right)  $ (continious functions vanishing at the
boundary) with $\psi_{t},\partial_{x}\psi\in L_{T}^{2}\left(  L_{loc}%
^{2}\right)  $. In view of the regularity properties of $X\left(  t,x\right)
$ it follows that for any $\psi\in\mathcal{D}\left(  (0,T)\times
\mathbb{R}\right)  $, $\psi\circ X^{-1}$ can be used as a test function. Using
this along with the fact that $X\left(  t,x\right)  $ is a homeomorphism for
all $t$, we write that%
\begin{align*}
0  & =\int_{0}^{T}\int_{\mathbb{R}}\rho(\partial_{t}\psi\circ X^{-1})+\rho
u\partial_{x}(\psi\circ X^{-1})dxdt\\
& =\int_{0}^{T}\int_{\mathbb{R}}\tilde{\rho}\left(  \widetilde{\partial
_{t}(\psi\circ X^{-1})}+\tilde{u}\widetilde{\partial_{x}(\psi\circ X^{-1}%
})\right)  \partial_{x}X\left(  t,x\right)  dxdt\\
& =\int_{0}^{T}\int_{\mathbb{R}}\tilde{\rho}\partial_{x}X\partial_{t}\psi
\end{align*}
witch translates into
\begin{equation}
\dfrac{d}{dt}\left(  \dfrac{\partial X}{\partial x}\tilde{\rho}\right)
=0.\label{rho_lag_equation}%
\end{equation}
Prooceding in a symilar manner, we get that the velocity's equation rewrites
as%
\begin{equation}
\rho_{0}\left(  x\right)  \partial_{t}\widetilde{u}-\partial_{x}\left(
\left(  \frac{\partial X}{\partial x}\right)  ^{-1}\mu(\widetilde{\rho
})\,\partial_{x}\widetilde{u}\right)  +\partial_{x}P(\widetilde{\rho
})=0.\label{u_lag_equation}%
\end{equation}
Putting together equations, $\left(  \text{\ref{rho_lag_equation}}\right)  $
and $\left(  \text{\ref{u_lag_equation}}\right)  $ we deduce that the system
$\left(  \text{\ref{Sistem_appendix}}\right)  $ can be writen in lagrangian
coordinates as:%
\begin{equation}
\left\{
\begin{array}
[c]{l}%
\dfrac{d}{dt}\left(  \dfrac{\partial X}{\partial x}\tilde{\rho}\right)  =0,\\
\rho_{0}\left(  x\right)  \partial_{t}\widetilde{u}-\partial_{x}\left(
\left(  \dfrac{\partial X}{\partial x}\right)  ^{-1}\mu(\widetilde{\rho
})\,\partial_{x}\,\widetilde{u}\right)  +\partial_{x}P\left(  \widetilde{\rho
}\right)  =0,\\
X\left(  t,x\right)  =x+%
{\displaystyle\int_{0}^{t}}
\widetilde{u}\left(  \tau,x\right)  d\tau,
\end{array}
\right.  \label{forma1}%
\end{equation}
or, equivalently%
\begin{equation}
\left\{
\begin{array}
[c]{l}%
\dfrac{d}{dt}\left(  \dfrac{\partial X}{\partial x}\widetilde{\rho}\right)
=0,\\
\rho_{0}\left(  x\right)  \partial_{t}\widetilde{u}-\partial_{x}\left(
\dfrac{\tilde{\rho}\mu(\widetilde{\rho})}{\rho_{0}}\,\partial_{x}%
\,\widetilde{u}\right)  +\partial_{x}P\left(  \widetilde{\rho}\right)  =0,\\
X\left(  t,x\right)  =x+%
{\displaystyle\int_{0}^{t}}
\widetilde{u}\left(  \tau,x\right)  d\tau,
\end{array}
\right.  \label{NS_lag}%
\end{equation}
Let us close this appendix observing that if we dispose of an inequality of
the following type (it is the case in our case, see (\ref{1.17})):
\begin{equation}
C(t)^{-1}\leq\widetilde{\rho}(t,x)\leq C\left(  t\right)  \label{dens}%
\end{equation}
then one may obtain from (\ref{NS_lag}) that%
\begin{equation}
C(t)^{-1}\inf\rho_{0}\leq\frac{\partial X}{\partial x}\left(  t,x\right)  \leq
C\left(  t\right)  \sup\rho_{0},\label{flow_ineg1}%
\end{equation}
along with
\begin{equation}
\frac{C(t)}{\inf\rho_{0}}\geq\left(  \frac{\partial X}{\partial x}\left(
t,x\right)  \right)  ^{-1}=\frac{\widetilde{\rho}(t,x)}{\rho_{0}(x)}\geq
\frac{C(t)^{-1}}{\sup\rho_{0}}.\label{flow_ineg2}%
\end{equation}

\subsection*{Sketch of the proof of the Theorem \ref{Cons}}
In this section, we are just giving a sketch of the proof of the blow-up criterion. The part concerning the existence of strong solution in finite time is classical.
We begin by observing that the Navier-Stokes system can be writen under the
following form:%
\begin{equation}
\left\{ 
\begin{array}{l}
\partial _{t}u+2u\partial _{x}u-\partial _{x}\left( \rho ^{\alpha
-1}\partial _{x}u\right) =v\partial _{x}u-\gamma \rho ^{\gamma -\alpha
}\left( v-u\right) , \\ 
\partial _{t}v+u\partial _{x}v=-\gamma \rho ^{\gamma -\alpha +1}\left(
v-u\right) 
\end{array}%
\right.   \label{sistem_Hs0}
\end{equation}%
Let us recall a classical product law in Sobolev spaces along with the
Kato-Ponce comutator estimate

\begin{lemme}[Kato-Ponce]
\label{lom}
The following estimates holds true for $s>0$ with ${\cal F}\Lambda_s f(\xi)=|\xi|^s{\cal F}f(\xi)$ for $f$ a temperated distribution:%
\begin{align}
\left\Vert \Lambda _{s}\left( fg\right) \right\Vert _{L^{2}}& \leq
\left\Vert f\right\Vert _{L^{\infty }}\left\Vert \Lambda _{s}g\right\Vert
_{L^{2}}+\left\Vert g\right\Vert _{L^{\infty }}\left\Vert \Lambda
_{s}f\right\Vert _{L^{2}} \\
\left\Vert \Lambda _{s}\left( f\partial _{x}g\right) -f\Lambda _{s}\partial
_{x}g\right\Vert _{L^{2}}& \leq C\left( \left\Vert \partial _{x}f\right\Vert
_{L^{\infty }}\left\Vert \Lambda _{s}g\right\Vert _{L^{2}}+\left\Vert
\Lambda _{s}f\right\Vert _{L^{2}}\left\Vert \partial _{x}g\right\Vert
_{L^{\infty }}\right) 
\end{align}
\end{lemme}
In the sequel we wish to describe how to preserve all along the time the $H^s$ norm of $u$ and $\rho-1$ for $s>\frac{3}{2}$.\\
We rewrite the system $\left( \text{\ref{sistem_Hs0}}\right) $ as%
\begin{equation}
\left\{ 
\begin{array}{l}
\partial _{t}\Lambda _{s}u+2u\partial _{x}\Lambda _{s}u-\partial _{x}\left(
\rho ^{\alpha -1}\partial _{x}\Lambda _{s}u\right) =\Lambda _{s}\left(
v\partial _{x}u\right) -\gamma \Lambda _{s}\left( \rho ^{\gamma -\alpha
}\left( v-u\right) \right)  \\ 
+2\left[ \Lambda _{s},u\right] \partial _{x}u+\partial _{x}\left( \left[
\rho ^{\alpha -1},\Lambda _{s}\right] \partial _{x}u\right) , \\ 
\partial _{t}\Lambda _{s}v+u\partial _{x}\Lambda _{s}v=-\gamma \Lambda
_{s}\left( \rho ^{\gamma -\alpha +1}\left( v-u\right) \right) +\left[
\Lambda _{s},u\right] \partial _{x}v%
\end{array}%
\right.   \label{sistem_Hs}
\end{equation}%
Multiply the first equation with $\Lambda _{s}u$ and integrate over $\R$, we get that: 
\begin{equation}
\begin{aligned}
&\frac{1}{2}\frac{d}{dt}\int_{\mathbb{R}}\left\vert \Lambda _{s}u\right\vert
^{2}+\int_{\mathbb{R}}\rho ^{\alpha -1}\left\vert \partial _{x}\Lambda
_{s}u\right\vert ^{2} =\int_{\mathbb{R}}\partial _{x}u\left\vert \Lambda
_{s}u\right\vert ^{2}+\int_{\mathbb{R}}\Lambda _{s}\left( v\partial
_{x}u\right) \Lambda _{s}u\\
&-\gamma \int_{\mathbb{R}}\Lambda _{s}\left( \rho
^{\gamma -\alpha }\left( v-u\right) \right) \Lambda _{s}u+\int_{\mathbb{R}}2\left[ \Lambda _{s},u\right] \partial _{x}u\Lambda
_{s}u+\int_{\mathbb{R}}\partial _{x}\left( \left[ \rho ^{\alpha -1},\Lambda
_{s}\right] \partial _{x}u\right) \Lambda _{s}u.  
\end{aligned}
\label{first_equation}
\end{equation}%
Multipliyng the second equation of $\left( \text{\ref{sistem_Hs}}%
\right) $ with $\Lambda _{s}v$ we obtain that: 
\begin{equation}
\frac{1}{2}\frac{d}{dt}\int_{\mathbb{R}}\left\vert \Lambda _{s}v\right\vert
^{2}=\frac{1}{2}\int_{\mathbb{R}}\partial _{x}u\left\vert \Lambda
_{s}v\right\vert ^{2}-\gamma \int_{\mathbb{R}}\Lambda _{s}\left( \rho
^{\gamma -\alpha +1}\left( v-u\right) \right) \Lambda _{s}v+\int_{\mathbb{R}}%
\left[ \Lambda _{s},u\right] \partial _{x}v\Lambda _{s}v.
\label{second_equation}
\end{equation}%
If we add up $\left( \text{\ref{first_equation}}\right) $ and $\left( \text{%
\ref{second_equation}}\right) $, it yields that%
\begin{eqnarray}
&&\frac{1}{2}\frac{d}{dt}\int_{\mathbb{R}}\left\{ \left\vert \Lambda
_{s}u\right\vert ^{2}+\left\vert \Lambda _{s}v\right\vert ^{2}\right\}
+\int_{\mathbb{R}}\rho ^{\alpha -1}\left\vert \partial _{x}\Lambda
_{s}u\right\vert ^{2}  \nonumber \\
&&=\int_{\mathbb{R}}\partial _{x}u\left\vert \Lambda _{s}u\right\vert
^{2}+\int_{\mathbb{R}}\Lambda _{s}\left( v\partial _{x}u\right) \Lambda
_{s}u-\gamma \int_{\mathbb{R}}\Lambda _{s}\left( \rho ^{\gamma -\alpha
}\left( v-u\right) \right) \Lambda _{s}u  \nonumber \\
&&+\int_{\mathbb{R}}2\left[ \Lambda _{s},u\right] \partial _{x}u\Lambda
_{s}u+\int_{\mathbb{R}}\partial _{x}\left( \left[ \rho ^{\alpha -1},\Lambda
_{s}\right] \partial _{x}u\right) \Lambda _{s}u  \nonumber \\
&&\frac{1}{2}\int_{\mathbb{R}}\partial _{x}u\left\vert \Lambda
_{s}v\right\vert ^{2}-\gamma \int_{\mathbb{R}}\Lambda _{s}\left( \rho
^{\gamma -\alpha +1}\left( v-u\right) \right) \Lambda _{s}v+\int_{\mathbb{R}}%
\left[ \Lambda _{s},u\right] \partial _{x}v\Lambda _{s}v  \label{bilan_semi}
\end{eqnarray}%
In the following lines, we analyse the different terms appearing in the left
hand side of $\left( \text{\ref{bilan_semi}}\right) $. The first two terms
are treated in the following manner using Lemma \ref{lom}:%
\begin{equation}
\begin{aligned}
&\int_{\mathbb{R}}\partial _{x}u\left\vert \Lambda _{s}u\right\vert
^{2}+\int_{\mathbb{R}}\Lambda _{s}\left( v\partial _{x}u\right) \Lambda _{s}u \\
&\lesssim \left\Vert \partial _{x}u\right\Vert _{L^{\infty }}\left\Vert
\Lambda _{s}u\right\Vert _{L^{2}}^{2}+\left\Vert \partial _{x}u\right\Vert
_{L^{\infty }}\left\Vert \Lambda _{s}v\right\Vert _{L^{2}} \left\Vert \Lambda
_{s}u\right\Vert _{L^{2}}+\left\Vert v\right\Vert _{L^{\infty }}\left\Vert
\partial _{x}\Lambda _{s}u\right\Vert _{L^{2}}\left\Vert \Lambda
_{s}u\right\Vert _{L^{2}}  \\
&\lesssim \left\Vert \partial _{x}u\right\Vert _{L^{\infty }}\left\Vert
\Lambda _{s}u\right\Vert _{L^{2}}^{2}+\left\Vert \partial _{x}u\right\Vert
_{L^{\infty }}\left\Vert \Lambda _{s}v\right\Vert  _{L^{2}}\left\Vert \Lambda
_{s}u\right\Vert _{L^{2}}+\left\Vert v\right\Vert _{L^{\infty }}\left\Vert
\rho ^{1-\alpha }\right\Vert _{L^{\infty }}^{\frac{1}{2}}\left\Vert \rho
^{\alpha -1}\partial _{x}\Lambda _{s}u\right\Vert _{L^{2}}\left\Vert \Lambda
_{s}u\right\Vert _{L^{2}}  \\
&\leq C\left\Vert \partial _{x}u\right\Vert _{L^{\infty }}\left\Vert
\Lambda _{s}u\right\Vert _{L^{2}}^{2}+C\left\Vert \partial _{x}u\right\Vert
_{L^{\infty }}\left\Vert \Lambda _{s}v\right\Vert _{L^{2}} \left\Vert \Lambda
_{s}u\right\Vert _{L^{2}}+C\left\Vert v\right\Vert _{L^{\infty
}}^{2}\left\Vert \rho ^{1-\alpha }\right\Vert _{L^{\infty }}\left\Vert
\Lambda _{s}u\right\Vert _{L^{2}}^{2}\\
&\hspace{10cm}+\frac{1}{8}\left\Vert \rho ^{\frac{%
\alpha -1}{2}}\partial _{x}\Lambda _{s}u\right\Vert _{L^{2}}^{2}.
\end{aligned}
\label{T1_T2}
\end{equation}%
The third term can be treated as follows: 
\begin{equation}
\begin{aligned}
&\int_{\mathbb{R}}\Lambda _{s}\left( \rho ^{\gamma -\alpha }\left(
v-u\right) \right) \Lambda _{s}u  \lesssim \big( \left\Vert \rho ^{\gamma -\alpha }\right\Vert _{L^{\infty
}}\left( \left\Vert \Lambda _{s}v\right\Vert _{L^{2}}+\left\Vert \Lambda
_{s}u\right\Vert _{L^{2}}\right) \left\Vert \Lambda _{s}u\right\Vert
_{L^{2}}\\
&\hspace{5cm}+\left \Vert \left( v-u\right) \right\Vert _{L^{\infty }}\left\Vert
\Lambda _{s}\left( \rho ^{\gamma -\alpha }-1\right) \right\Vert
_{L^{2}}\big) \left\Vert \Lambda _{s}u\right\Vert _{L^{2}}. 
\end{aligned}
 \label{T3}
\end{equation}%
We have for the fourth term using Lemma \ref{lom}:%
\begin{equation}
\int_{\mathbb{R}}2\left[ \Lambda _{s},u\right] \partial _{x}u\Lambda
_{s}u\lesssim \left\Vert \partial _{x}u\right\Vert _{L^{\infty }}\left\Vert
\Lambda _{s}u\right\Vert _{L^{2}}^{2}  \label{T4}
\end{equation}%
The fifth term :%
\begin{equation}
\begin{aligned}
&\int_{\mathbb{R}}\partial _{x}\left( \left[ \rho ^{\alpha -1},\Lambda _{s}%
\right] \partial _{x}u\right) \Lambda _{s}u   \\
&\leq C\left\Vert \left[ \rho ^{\alpha -1},\Lambda _{s}\right] \partial
_{x}u\right\Vert _{L^{2}}\left\Vert \partial _{x}\Lambda _{s}u\right\Vert
_{L^{2}}   \\
&\leq C\left\Vert \rho ^{\alpha -1}\right\Vert _{L^{\infty }}\left(
\left\Vert \partial _{x}\rho ^{\alpha -1}\right\Vert _{L^{\infty
}}\left\Vert \Lambda _{s}u\right\Vert _{L^{2}}+\left\Vert \partial
_{x}u\right\Vert _{L^{\infty }}\left\Vert \Lambda _{s}\left( \rho ^{\alpha
-1}-1\right) \right\Vert _{L^{2}}\right) ^{2}+\frac{1}{8}\left\Vert \rho ^{%
\frac{\alpha -1}{2}}\partial _{x}\Lambda _{s}u\right\Vert _{L^{2}}^{2} 
 \\
&\leq C\left\Vert \rho ^{\alpha -1}\right\Vert _{L^{\infty }}\left(
\left\Vert v-u\right\Vert _{L^{\infty }}\left\Vert \Lambda _{s}u\right\Vert
_{L^{2}}+\left\Vert \partial _{x}u\right\Vert _{L^{\infty }}\left\Vert
\Lambda _{s}\left( \rho ^{\alpha -1}-1\right) \right\Vert _{L^{2}}\right)
^{2}+\frac{1}{8}\left\Vert \rho ^{\frac{\alpha -1}{2}}\partial _{x}\Lambda
_{s}u\right\Vert _{L^{2}}^{2} .
\end{aligned}
\label{T5}
\end{equation}
We skip the sixth term. Seventh term :%
\begin{equation}
\begin{aligned}
&\int_{\mathbb{R}}\Lambda _{s}\left( \rho ^{\gamma -\alpha +1}\left(
v-u\right) \right) \Lambda _{s}v  \lesssim \left\Vert \Lambda _{s}v\right\Vert _{L^{2}}\left\Vert \Lambda
_{s}u\right\Vert _{L^{2}}+\left\Vert \Lambda _{s}v\right\Vert _{L^{2}}^{2} 
 \\
&+\left( \left\Vert \rho ^{\gamma -\alpha +1}-1\right\Vert _{L^{\infty
}}\left\Vert \Lambda _{s}\left( v-u\right) \right\Vert _{L^{2}}+\left\Vert
v-u\right\Vert _{L^{\infty }}(\left\Vert \Lambda _{s}\left( \rho ^{\gamma
-\alpha +1}-1\right) \right\Vert _{L^{2}})\right) \left\Vert \Lambda
_{s}v\right\Vert _{L^{2}}
\end{aligned}  \label{T7}
\end{equation}%
Last term : 
\begin{equation}
\int_{\mathbb{R}}\left[ \Lambda _{s},u\right] \partial _{x}v\Lambda
_{s}v\lesssim \left( \left\Vert \partial _{x}u\right\Vert _{L^{\infty
}}\left\Vert \Lambda _{s}v\right\Vert _{L^{2}}+\left\Vert \partial
_{x}v\right\Vert _{L^{\infty }}\left\Vert \Lambda _{s}u\right\Vert
_{L^{2}}\right) \left\Vert \Lambda _{s}v\right\Vert _{L^{\acute{e}}}
\label{T8}
\end{equation}
Let us observe that in the estimates $\left( \text{\ref{T3}}\right) $ $%
\left( \text{\ref{T5}}\right) $ and $\left( \text{\ref{T7}}\right) $ we have
to treat the $H^{s}$-norm of $\rho ^{\gamma -\alpha }$,$\rho ^{\alpha -1}$
and $\rho ^{\gamma -\alpha +1}$ respectively. This is the objective of the
following lines. For each $\beta $, we may write that%
\[
\partial _{t}\rho ^{\beta }+u\partial _{x}\rho ^{\beta }=-\beta \rho ^{\beta
}\partial _{x}u.
\]%
and consequently 
\[
\partial _{t}\Lambda _{s}\left( \rho ^{\beta }-1\right) +u\partial
_{x}\Lambda _{s}\left( \rho ^{\beta }-1\right) =-\beta \Lambda _{s}\left(
\rho ^{\beta }\partial _{x}u\right) +\left[ \Lambda _{s},u\right] \partial
_{x}\left( \rho ^{\beta }-1\right) .
\]%
We get that%
\begin{eqnarray}
&&\frac{1}{2}\frac{d}{dt}\int_{\mathbb{R}}\left\vert \Lambda _{s}\left( \rho
^{\beta }-1\right) \right\vert ^{2} \leq (\frac{1}{2}+\beta )\left\Vert
\partial _{x}u\right\Vert _{L^{\infty }}\left\Vert \Lambda _{s}\left( \rho
^{\beta }-1\right) \right\Vert _{L^{2}}^{2}  \nonumber \\
&&+\beta \left( \left\Vert \partial _{x}u\right\Vert _{L^{\infty
}}\left\Vert \Lambda _{s}\left( \rho ^{\beta }-1\right) \right\Vert
_{L^{2}}^{2}+\left\Vert \left( \rho ^{\beta }-1\right) \right\Vert
_{L^{\infty }}\left\Vert \Lambda _{s}\partial _{x}u\right\Vert
_{L^{2}}\right) \left\Vert \Lambda _{s}\left( \rho ^{\beta }-1\right)
\right\Vert _{L^{2}}  \nonumber \\
&&+\left( \left\Vert \partial _{x}u\right\Vert _{L^{\infty }}\left\Vert
\Lambda _{s}\left( \rho ^{\beta }-1\right) \right\Vert _{L^{2}}+\left\Vert
\partial _{x}\rho ^{\beta }\right\Vert _{L^{\infty }}\left\Vert \Lambda
_{s}u\right\Vert _{L^{2}}\right) \left\Vert \Lambda _{s}\left( \rho ^{\beta
}-1\right) \right\Vert _{L^{2}}  \nonumber \\
&&\leq C_{\varepsilon }\left( \left\Vert \partial _{x}u\right\Vert
_{L^{\infty }}+\left\Vert \rho ^{1-\alpha }\right\Vert _{L^{\infty
}}^{2}\left\Vert \left( \rho ^{\beta }-1\right) \right\Vert _{L^{\infty
}}^{2}+\left\Vert v-u\right\Vert _{L^{\infty }}^{2}\right) \left\Vert
\Lambda _{s}\left( \rho ^{\beta }-1\right) \right\Vert _{L^{2}}^{2} 
\nonumber \\
&&+\left\Vert \Lambda _{s}u\right\Vert _{L^{2}}^{2}+\varepsilon \left\Vert
\rho ^{\frac{\alpha -1}{2}}\partial _{x}\Lambda _{s}u\right\Vert _{L^{2}}^{2}
\label{density_Hs}
\end{eqnarray}%
Thus putting togheter the estimates $\left( \text{\ref{T1_T2}}\right) $,$%
\left( \text{\ref{T3}}\right) $,$\left( \text{\ref{T4}}\right) $,$\left( 
\text{\ref{T5}}\right) $,$\left( \text{\ref{T7}}\right) $,$\left( \text{\ref%
{T8}}\right) $ and $\left( \text{\ref{density_Hs}}\right) $ for $\beta
=\alpha -1,$ $\gamma -\alpha ,$ $\gamma -\alpha +1$ we get that 
\begin{equation}
\begin{aligned}
&\int_{\mathbb{R}}\left\{ \left\vert \Lambda _{s}u\right\vert
^{2}+\left\vert \Lambda _{s}v\right\vert ^{2}+\left\vert \Lambda _{s}\left(
\rho ^{\alpha -1}-1\right) \right\vert ^{2}+\left\vert \Lambda _{s}\left(
\rho ^{\gamma -\alpha }-1\right) \right\vert ^{2}+\left\vert \Lambda
_{s}\left( \rho ^{\gamma -\alpha +1}-1\right) \right\vert ^{2}\right\}(t,x) dx
\\
&+\int_{0}^{t}\int_{\mathbb{R}}\rho ^{\alpha -1}\left\vert \partial
_{x}\Lambda _{s}u\right\vert ^{2}  (s,x) ds dx \\
&\leq C\left( u_{0},\rho _{0}\right) \exp \left( \int_{0}^{t}\left(
1+\left\Vert \left( \rho ,\frac{1}{\rho }\right) \right\Vert _{L^{\infty
}}\right) ^{\delta }\left( 1+\left\Vert \left( u,v,\partial _{x}u,\partial
_{x}v\right) \right\Vert _{L^{\infty }}^{2}\right) \right) 
\end{aligned}
\label{Basic_explosion}
\end{equation}%
with $\delta $ depending on $\alpha $ and $\gamma $. We mention also that $C\left( u_{0},\rho _{0}\right)$ depends on $\|u_0\|_{H^s}$, $\|\rho_0-1\|_{H^s}$, $\|\rho_0\|_{L^\infty}$ and $\|\frac{1}{\rho_0}\|_{L^\infty}$. Next, let us analyse in more detail the equation of $v$ : 
\[
\partial _{t}v+u\partial _{x}v=-\gamma \rho ^{\gamma -\alpha +1}\left(
v-u\right) 
\]%
We get that%
\[
\left\Vert v\right\Vert _{L_{t}^{\infty }(L_{x}^{\infty })}\leq \left\Vert
v_{0}\right\Vert _{L_{x}^{\infty }}+\left\Vert \rho ^{\gamma -\alpha
+1}\right\Vert _{L_{t}^{\infty }(L_{x}^{\infty })}\int_{0}^{t}\left\Vert
u(s,\cdot)\right\Vert _{L^{\infty }}ds.
\]%
Moreover, writing the equation of $\partial _{x}v$ we see that%
\[
\left( \partial _{x}v\right) _{t}+u\partial _{x}\left( \partial _{x}v\right)
+\left( \partial _{x}u+\gamma \rho ^{\gamma -\alpha }\right) \partial
_{x}v=\gamma \rho ^{\gamma -\alpha }\partial _{x}u-\gamma \left( \gamma
-\alpha \right) \rho ^{\gamma -2\alpha +2}\left( v-u\right) ^{2}
\]%
From which we deduce that%
\[
\left\Vert \partial _{x}v\right\Vert _{L_{t}^{\infty }(L_{x}^{\infty
})}\lesssim \left\Vert \partial _{x}v_{0}\right\Vert _{L_{x}^{\infty }}+\psi
\left( \left( 1+\left\Vert \left( \rho ,\frac{1}{\rho }\right) \right\Vert
_{L_{t}^{\infty }(L^{\infty })}\right) ^{\delta _{1}}\left(
1+\int_{0}^{t}\left\Vert \left( u,\partial _{x}u\right) \right\Vert
_{L^{\infty }}\right) \right) 
\]%
with $\psi \left( r\right) =r\exp r$ and $\delta _{1}$ depending on $\gamma $
and $\alpha $. Moreover, the Bresch-Desjardins entropy allows a uniform control on $%
\left\Vert \rho \right\Vert _{L_{t}^{\infty }(L_{x}^{\infty })}$.
Let us denote by: 
\[
\tilde{A}\left( \rho ,u\right) \left( t\right) =\int_{\mathbb{R}}\rho
^{\alpha }\left( t\right) \left( \partial _{x}u\right) ^{2}\left( t\right)
+\int_{0}^{t}\int_{\mathbb{R}}\rho \dot{u}^{2}.
\]%
Using the same techniques as in the section on the Hoff estimates, we may show that%
\[
\tilde{A}\left( \rho ,u\right) \left( t\right) \leq C_{0}\exp \left( t\left(
1+\left\Vert \frac{1}{\rho }\right\Vert _{L_{t}^{\infty }(L^{\infty
})}\right) ^{\delta _{2}}\right) 
\]%
which, in turn, ensures a control on $\left\Vert \partial _{x}u\right\Vert
_{L^{2}_t(L^{\infty })}$ provided that we control $\|\frac{1}{\rho}\|_{L^\infty_t(L^\infty)}$. To summarize:
\begin{itemize}
\item The Bresch-Desjardins entropy provides control on $\left\Vert
\rho \right\Vert _{L_{t}^{\infty }(L_{x}^{\infty })}$ for any $t>0$, 

\item $\left\Vert \left( v,\partial _{x}v\right) \right\Vert _{L_t^{\infty
}(L^{\infty })}$ is controlled by $\left\Vert \left( u,\partial _{x}u\right)
\right\Vert _{L^{1}\left( L^{\infty }\right) }$ and $\left\Vert \left(
\rho ,\frac{1}{\rho }\right) \right\Vert _{L_{t}^{\infty }(L^{\infty })},$

\item The Hoff-type estimates ensure that $\left\Vert \partial
_{x}u\right\Vert _{L_{t}^{2}(L_{x}^{\infty })}$ is  controlled by $%
\left\Vert \left( \rho ,\frac{1}{\rho }\right) \right\Vert _{L_{t}^{\infty
}(L^{\infty })}$

\item Using the basic energy estimate we obtain that $\left\Vert u\right\Vert
_{L_{t}^{2}(L_{x}^{\infty })}$ is  controlled by $\left\Vert \left(
\rho ,\frac{1}{\rho }\right) \right\Vert _{L_{t}^{\infty }(L^{\infty })}.$
\end{itemize}
Taking into account the estimate $\left( \text{\ref{Basic_explosion}}\right) 
$ we get that for any $T>0$ and any $s>\frac{3}{2}$ the $H^{s}$-Sobolev norm of $%
\left( u,v,\rho-1 \right) $ is uniformly controlled by $\left\Vert \frac{1}{%
\rho }\right\Vert _{L_{t}^{\infty }(L^{\infty })}.$

\section*{Acknowledgements}

CB has been partially funded by the ANR project SingFlows ANR-18-CE40-0027-01.
BH has been partially funded by the ANR project INFAMIE ANR-15-CE40-0011. This
work was realized during the secondment of BH in the ANGE Inria team.

\end{document}